\documentclass{amsart}

\newtheorem{thm}{Theorem}

\newtheorem{lem}[thm]{Lemma}
\newtheorem{cor}[thm]{Corollary}

\newtheorem{prop}[thm]{Proposition}

\theoremstyle{definition}
\newtheorem{defn}[thm]{Definition}

\newtheorem{exmp}[thm]{Example}

\newtheorem{prob}[thm]{Problem}

\newtheorem{rem}[thm]{Remark}
\newtheorem{ack}{Acknowledgments}

\newtheorem{defn-thm}[thm]{Definition--Theorem}  

\theoremstyle{remark}

\newcommand{\CX}{\mbox{${\calx}\hspace{-.8em}-\,$}}

\newcommand{\hl}{\\\hline}

\def\fract#1#2{\raise4pt\hbox{$ #1 \atop #2 $}}
\def\decdnar#1{\phantom{\hbox{$\scriptstyle{#1}$}}
\left\downarrow\vbox{\vskip15pt\hbox{$\scriptstyle{#1}$}}\right.}

\def\bfa{{\bf a}}

\def\bfw{{\bf w}}
\def\bfx{{\bf x}}
\def\bfy{{\bf y}}
\def\bfz{{\bf z}}

\def\calc{{\mathcal C}}
\def\calo{{\mathcal O}}

\def\cald{{\mathcal D}}
\def\cale{{\mathcal E}}
\def\calf{{\mathcal F}}

\def\call{{\mathcal L}}

\def\calo{{\mathcal O}}

\def\cals{{\mathcal S}}

\def\calw{{\mathcal W}}

\def\calz{{\mathcal Z}}
\def\calx{{\mathcal X}}

\def\bba{{\mathbb A}}

\def\bbc{{\mathbb C}}
\def\bbd{{\mathbb D}}

\def\bbh{{\mathbb H}}
\def\bbi{{\mathbb I}}

\def\bbo{{\mathbb O}}
\def\bbp{{\mathbb P}}
\def\bbq{{\mathbb Q}}
\def\bbr{{\mathbb R}}

\def\bbt{{\mathbb T}}

\def\bbz{{\mathbb Z}}

\def\gra{\alpha}
\def\grb{\beta}

\def\grd{\delta}
\def\gre{\epsilon}

\def\gri{\iota}

\def\grl{\lambda}

\def\gro{\omega}

\def\grr{\rho}
\def\grs{\sigma}

\def\grz{\zeta}

\def\grG{\Gamma}

\def\grL{\Lambda}
\def\grO{\Omega}

\def\grS{\Sigma}

\def\gsp1{{\mathfrak s}{\mathfrak p}(1)}

\def\teta{\tilde{\eta}}

\def\ker{\hbox{ker}}

\def\Ric{\hbox{Ric}}
\def\barj{\bar{j}}

\def\ga{{\mathfrak a}}

\def\gt{{\mathfrak t}}
\def\gu{{\mathfrak u}}

\def\gF{{\mathfrak F}}

\def\gI{{\mathfrak I}}

\def\gT{{\mathfrak T}}

\def\tg{\tilde{{g}}}

\def\txi{\tilde{{\xi}}}
\def\teta{\tilde{{\eta}}}
\def\tPhi{\tilde{{\Phi}}}
\def\tg{\tilde{g}}
\def\tM{\tilde{M}}

\def\teta{\tilde{\eta}}
\def\ker{\hbox{ker}}

\def\Ric{{\rm Ric}}

\def\semidirright{\times\kern-8.0pt\hbox{$\triangleleft$}}
\def\hsmash{\triangleright\kern-4.4pt\raise.66pt\hbox{$\scriptstyle{<}$}}

\def\la#1{\hbox to #1pc{\leftarrowfill}}
\def\ra#1{\hbox to #1pc{\rightarrowfill}}
\def\Se{Sasakian-Einstein }

\begin{document}
\bibliographystyle{amsalpha}

\title{On Eta-Einstein Sasakian Geometry}
\author{Charles P. Boyer, Krzysztof Galicki and Paola Matzeu}
\address{CPB: Department of Mathematics \& Statistics,
University of New Mexico, Albuquerque, NM 87131.}
\email{cboyer@math.unm.edu}
\address{KG: Max-Planck-Institut f\"ur Mathematik,
D53111 Bonn, Germany; on leave from
Department of Mathematics and Statistics,
University of New Mexico,
Albuquerque, NM 87131}
\email{galicki@math.unm.edu}
\address{PM: Dipartimento di Matematica e Informatica,
Universit\'a di Cagliari, Cagliari, Italy.}
\email{matzeu@vaxca1.unica.it}
\begin{abstract}
A compact quasi-regular Sasakian manifold $M$ is foliated by one-dimensional leaves
and the transverse space of this characteristic foliation is necessarily
a compact K\"ahler orbifold $\calz$. In the case when the transverse space $\calz$ is also
Einstein the corresponding Sasakian manifold $M$ is said to be Sasakian $\eta$-Einstein. In
this article
we study $\eta$-Einstein geometry as a class of distinguished Riemannian metrics
on contact metric manifolds. In particular, we use a previous solution of the Calabi
problem in the context of Sasakian geometry to prove the existence of $\eta$-Einstein
structures on many
different compact manifolds, including exotic spheres. We also relate these results to the
existence of Einstein-Weyl and Lorenzian Sasakian-Einstein structures.
\end{abstract}
\subjclass{Primary 53C25}
\keywords{Sasakian geometry, Einstein-Weyl geometry, K\"ahler-Einstein metrics, Calabi
problem, orbifolds, contact structures}

\maketitle

\bigskip
\section{Introduction}
\medskip
The purpose of this paper is to study a special kind of Riemannian
metric on Sasakian manifolds. A Sasakian manifold $M$ of dimension
$2n+1$ with a Sasakian structure $\cals=(\xi,\eta,\Phi,g)$ is said
to be $\eta$-Einstein if the Ricci curvature tensor of the metric
$g$ satisfies the equation ${\rm Ric}_g=\lambda
g+\nu\eta\otimes\eta$ for some constants $\lambda,\nu\in\bbr$.
These metrics were introduced and studied by Okumura \cite{Oku62}
and then named by Sasaki \cite{Sas65} in his lecture notes in
1965. Okumura assumed that both $\grl$ and $\nu$ are functions,
and then proved, similar to the case of Einstein metrics, that
they must be constant when $n>1.$ Obviously $\nu=0$ reduces to the
more familiar Sasakian-Einstein condition. In general,
$\lambda+\nu=2n$ and every Sasakian $\eta$-Einstein manifold is
necessarily of constant scalar curvature $s=2n(1+\lambda)$.

It is well-known that Sasakian-Einstein metrics are necessarily
positive with Einstein constant equal to ${\rm dim}(M)-1$. On the
other hand, on K\"ahler manifolds one naturally considers
K\"ahler-Einstein metrics with Einstein constant of any sign. As
every Sasakian manifold comes with a characteristic foliation
$\calf_\xi$ whose transverse geometry is K\"ahler, the notion of
the basic first Chern class $c_1^B$ leads one to define null,
positive, or negative Sasakian structures in analogy with K\"ahler
geometry. One simply asks that $c_1^B$ either vanish or be
represented by a basic positive or negative $(1,1)$-form,
respectively. In such a context Sasakian $\eta$-Einstein metrics
are the natural analogues of K\"ahler-Einstein metrics.

In the $c_1^B>0$ case the transverse geometry is Fano and the
$\eta$-Einstein condition implies further that the transverse
geometry is K\"ahler-Einstein of positive scalar curvature. It was
first observed by Tanno that a positive Sasakian manifold $M$ with
an $\eta$-Einstein metric $g$ admits another metric $g'$ which is
Sasakian-Einstein \cite{Tan79}. The respective Sasakian structures
$\cals$ and $\cals'$ are related by the so-called `transverse' or
`$\cald$-homothety' transformation (see Section 3 for precise
definition). Tanno used this simple observation to show that a
unit tangent bundle $M=T_1S^n$ of any $n$-sphere admits a
homogeneous Sasakian-Einstein structure. Actually, transverse
homotheties define a family of Sasakian $\eta$-Einstein structures
for each scale $a\in\bbr^+$. This family can be represented by the
unique Sasakian-Einstein structure which then can be either
``squashed" or "stretched". It turns out that the squashed
Sasakian $\eta$-Einstein structures can be used to define an
associated Einstein-Weyl geometry in a canonical way \cite{Nar93,
Nar97, Nar98, Hig93, PeSw93}. For instance, in 3 dimensions the
only compact simply connected Sasakian-Einstein manifold is the
round $S^3$ and, hence, the squashed $\eta$-Einstein metric is
simply the Berger metric. The latter is a well-known example of a
non-trivial Einstein-Weyl manifold. It follows that
Sasakian-Einstein geometry offers quite a powerful tool in
constructing odd-dimensional examples of Einstein-Weyl manifolds.
In particular, many compact spin 5-manifolds admit families of
Einstein-Weyl structures. Similarly, both the standard as well as
some exotic spheres in odd dimensions also admit Einstein-Weyl
structures.

The $c_1^B=0$ case is an example of the so-called transverse
Calabi-Yau structure. As a consequence of a transverse version
\cite{ElK} of Yau's famous theorem, every null Sasakian structure
can be deformed to a Sasakian $\eta$-Einstein structure which is
transverse Calabi-Yau. We exhibit interesting examples of such
structures on certain compact Sasakian manifolds. Interestingly,
some of these examples are also Einstein-Weyl as observed by
Narita in \cite{Nar93}.

In the $c_1^B<0$ case an orbifold version of the theorem of Aubin
and Yau shows that every negative Sasakian structure $\cals$ can
be deformed to a Sasakian $\eta$-Einstein structure. As a
consequence, negative Sasakian $\eta$-Einstein metrics can be
easily found on many compact contact spin manifolds. Remarkably,
it follows that (possibly all) homotopy spheres which bound
parallelizable manifolds admit infinitely many inequivalent
negative Sasakian $\eta$-Einstein structures. Hence, for example,
$S^5$ admits both positive and negative Sasakian $\eta$-Einstein
structures, but no null Sasakian structures. Both positive and
negative Sasakian $\eta$-Einstein structures exist on many compact
simply connected spin  5-manifolds, including $\#k(S^2\times S^3)$
for certain $k$ (positive Sasakian structures exist for all $k$)
and certain rational homology 5-spheres. Some of these manifolds,
for example $\#20(S^2\times S^3)$ admit all three types of
Sasakian $\eta$-Einstein structures. One additional interesting
feature of the $c_1^B<0$ case is that every negative Sasakian
$\eta$-Einstein structure yields a family of Lorentzian Sasakian
$\eta$-Einstein metrics \cite{Bau00, Boh03}. The metrics in the
family are related by a transverse homothety and, like in the
positive case, there is a unique Lorentzian Sasakian-Einstein
metric in the family with Einstein constant $1-{\rm dim}(M)$.
Hence, for instance, $S^5$ admits infinitely many Lorentzian
Sasakian-Einstein structures.

By virtue of admitting real Killing spinors, Sasakian-Einstein
manifolds have recently received a lot of attention in physics.
The physicists' interest in five and seven-dimensional manifolds
admitting real Killing spinors dates back to the early eighties,
when Kaluza-Klein models played a central role in the supergravity
theory. Today's renewed interest in these manifolds has to do with
the so-called $p$-brane solutions in superstring theory. These
$p$-branes, ``near the horizon" are modeled by the
pseudo-Riemannian geometry of the product ${\rm adS}_{p+2}\times
M$, where ${\rm adS}_{p+2}$ is the $(p+2)$-dimensional
anti-de-Sitter space (a Lorentzian version of a space of constant
sectional curvature) and $(M, g)$ is a Riemannian manifold of
dimension $d=D-p-2$. Here $D$ is the dimension of the original
supersymmetric theory. In the most interesting cases of M2-branes,
M5-branes, and D3-branes $D$ equals to either 11 (M$p$-branes of
M-theory) or 10 (D$p$-branes in type IIA or type IIB string
theory). Any residual supersymmetry forces $M$ to admit real
Killing spinors. For example, the case of D3-branes of string
theory the relevant near horizon geometry is that of ${\rm
adS}_{5}\times M$, where $M$ is a \Se 5-manifold. The D3-brane
solution interpolates between ${\rm adS}_{5}\times M$ and
$\bbr^{3,1}\times \calc(M)$, where the cone $\calc(M)$ is a
Calabi-Yau threefold. In its original version the Maladacena
conjecture (also known as AdS/CFT duality) states that the 't
Hooft large $n$ limit of $N=4$ supersymmetric Yang-Mills theory
with gauge group $SU(n)$ is dual to type IIB superstring theory on
${\rm adS}_{5}\times S^5$ \cite{Mal99}. This conjecture has now
been examined for many known \Se metric in dimension 5 and 7. All
this has led to a remarkable discovery of new cohomogeneity one
\Se metrics \cite{GMSW04a, GMSW04b}. In dimension 5 the new
examples produce infinitely many toric \Se metric on $S^2\times
S^3$. The homogeneous \Se metric on $S^2\times S^3$ examined in
the context of the AdS/CFT duality by Klebanov and Witten
\cite{KlWi99} emerges as a special case of this construction
\cite{MaSp05, MaSpYau05}. These metrics are of equal if not greater
interest from the geometric point of view. They provide the first
examples of compact \Se manifolds which are not quasi-regular,
i.e., the space of leaves of the associated characteristic foliation
is not a topologically nicely behaved object such as an orbifold. As a consequence they are
also counterexamples to a conjecture made in 1994 by Cheeger and Tian
\cite{ChTia94}.

Deformations of Sasakian-Einstein metrics naturally lead to
Einstein-Weyl structures are important in general relativity
theory. Null and negative Sasakian $\eta$-Einstein manifolds have
not perhaps received as much attention. Yet, both in physics and
in mathematics, they often emerge in a somewhat implicit fashion.
For example, null Sasakian manifolds are orbifold bundles over
compact Calabi-Yau orbifolds. They are very natural
odd-dimensional companions of the Calabi-Yau spaces. Hence, for
example, 7-dimensional null Sasakian $\eta$-Einstein manifolds
should be important in the study of Calabi-Yau 3-folds with cyclic
orbifold singularities. In particular all problems regarding
famous Mirror Symmetry should have a translation in the language
of null Sasakian $\eta$-Einstein 7-manifolds.

Our article is organized as follows: In Section 2 we recall some
basic properties of the transverse geometry including invariants
of the characteristic foliation. Section 3 follows with
fundamentals about certain deformations of Sasakian structures. In
section 4 we begin the study of $\eta$-Einstein metrics and their
relation to the Sasakian Calabi problem. Here we also describe the
non-spin obstruction to the existence of $\eta$-Einstein metrics.
The following section describes some general properties of
Sasakian $\eta$-Einstein manifolds, in particular, we describe the
relation between negative Sasakian $\eta$-Einstein structures and
Lorentzian \Se structures. Many of the most interesting examples
of positive, null, and negative Sasakian structures can be found
on links of isolated hypersurface singularities which are
discussed in Section 6. In Sections 7 and 8 we use this fact to
introduce a large set of examples of Sasakian $\eta$-Einstein
structures. Finally, in the last section we discuss relations
between Einstein-Weyl and Sasakian $\eta$-Einstein structures.

\begin{ack} We thank Gang Tian for answering several questions.
We would also like to thank J\'anos Koll\'ar for numerous
discussions and for explaining
some of his recent results \cite{Kol04b} to us, and to Liviu Ornea for pointing out reference
\cite{GaOr98}.
CPB and KG were partially supported by the NSF under grant number
DMS-0203219. The authors would also like to thank Universit\`a di
Roma ``La Sapienza" for partial support where discussions on this
work initiated. The second named author would like to thank
Max-Planck-Institute in Bonn for hospitality and support. Much of
the work on this article was done during his stay there.
\end{ack}

\bigskip
\section{The Characteristic Foliation and Transverse Invariants}
\medskip
Recall that an {\it almost contact structure} on a differentiable
manifold $M$ is a triple $(\xi,\eta,\Phi),$ where $\Phi$ is a
tensor field of type $(1,1)$ (i.e., an endomorphism of $TM$),
$\xi$ is a vector field, and $\eta$ is a 1-form which satisfy
$$\eta(\xi)=1 ~~~\hbox{and}~~~ \Phi\circ \Phi= -\bbi + \xi\otimes\eta,$$
where $\bbi$ is the identity endomorphism on $TM.$ A smooth
manifold with such a structure is called an {\it almost contact
manifold}. A Riemannian metric $g$ on $M$ is said to be {\it
compatible} with the almost contact structure if for any vector
fields $X,Y$ on $M$ we have
$$g(\Phi X,\Phi Y)= g(X,Y)- \eta(X)\eta(Y).$$
\noindent An almost contact structure with a compatible metric is
called an {\it almost contact metric structure}. In case $\eta$ is
a contact form $(\xi,\eta,\Phi,g)$ is said to be a {\it contact
metric structure} on $M$.

A contact metric structure $(\xi,\eta,\Phi,g)$ is called {\it
K-contact} if $\xi$ is a Killing vector field of $g$ and it is
called {\it  Sasakian} if the metric cone
$(C(M),dr^2+r^2g,d(r^2\eta))$ is K\"ahler.

In this section we study the transverse geometry of the Riemannian
foliation $\calf_\xi$ defined on $M$ by the characteristic, or as
it is called {\bf Reeb} vector field $\xi.$ The transverse geometry of Sasakian structures
were discussed in part 3 of Sasaki's notes \cite{Sas68} in the regular case. This predated
the modern development of foliation theory which is the proper setting for the study of
Sasakian geometry.  We first make note of
some well-known properties. The foliation $\calf_\xi$ is one
dimensional whose leaves are geodesics with respect to the
Sasakian metric $g,$ and this metric is bundle-like.

Let $(M,\xi,\eta,\Phi,g)$ be a Sasakian manifold, and consider the
contact subbundle $\cald=\ker~\eta.$ There is an orthogonal
splitting of the tangent bundle as
\begin{equation}
TM=\cald \oplus L_\xi,
\end{equation}
where $L_\xi$ is the trivial line bundle generated by the Reeb
vector field $\xi.$ The contact subbundle $\cald$ is just the
normal bundle to the characteristic foliation $\calf_\xi$
generated by $\xi.$ It is naturally endowed with both a complex
structure $J=\Phi|\cald$ and a symplectic structure $d\eta.$
Hence, $(\cald,J,d\eta)$ gives $M$ a {\it transverse K\"ahler}
structure with K\"ahler form $d\eta$ and metric $g_\cald$ defined
by
\begin{equation}
g_\cald(X,Y)=d\eta(X,JY)
\end{equation}
which is related to the Sasakian metric $g$ by
\begin{equation}
g=g_\cald \oplus \eta\otimes \eta.
\end{equation}
Recall, see e.g. \cite{Ton}, that a smooth p-form $\gra$ on $M$ is
called {\it basic} if
\begin{equation}
\xi\rfloor \gra=0, \qquad \call_\xi\gra=0,
\end{equation}
and we let $\grL^p_B$ denote the sheaf of germs of basic p-forms
on $M,$ and by $\grO_B^p$ the set of global sections of $\grL^p_B$
on $M.$ The sheaf $\grL^p_B$ is a module under the ring,
$\grL^0_B,$ of germs of smooth basic functions on $M.$ We let
$C^\infty_B(M)=\grO^0_B$ denote global sections of $\grL^0_B,$
i.e. the ring of smooth basic functions on $M.$  Since exterior
differentiation preserves basic forms we get a de Rham complex
\begin{equation}
\cdots\ra{2.5}\grO_B^p\fract{d}{\ra{2.5}}\grO_B^{p+1}\ra{2.5}\cdots
\end{equation}
whose cohomology $H^*_B(\calf_\xi)$ is called the {\it basic
cohomology} of $(M,\calf_\xi).$ The basic cohomology ring
$H^*_B(\calf_\xi)$ is an invariant of the foliation $\calf_\xi$
and hence, of the Sasakian structure on $M.$ When $M$ is compact
it is related to the ordinary de Rham cohomology $H^*(M,\bbr)$ by
the long exact sequence, (see e.g. \cite{Ton})
\begin{equation}\label{exact}
\cdots\ra{2.5}H_B^p(\calf_\xi)\ra{2.5}H^p(M,\bbr)\fract{j_p}{\ra{2.5}}
H_B^{p-1}(\calf_\xi) \fract{\grd}{\ra{2.5}}
H^{p+1}_B(\calf_\xi)\ra{2.5}\cdots
\end{equation}
where $\grd$ is the connecting homomorphism given by
$\grd[\gra]=[d\eta\wedge \gra]=[d\eta]\cup[\gra],$ and $j_p$ is
the composition of the map induced by $\xi\rfloor$ with the well
known isomorphism $H^r(M,\bbr)\approx H^r(M,\bbr)^{\gT}$ where $\gT,$ a torus, is the
closure of the Reeb flow, and
$H^r(M,\bbr)^{\gT}$ is the $\gT$-invariant cohomology defined from
the  $\gT$-invariant r-forms $\grO^r(M)^{\gT}.$ We also note that
$d\eta$ is basic even though $\eta$ is not.

Next we exploit the fact that the transverse geometry is K\"ahler.
Let $\cald_\bbc$ denote the complexification of $\cald,$ and
decompose it into its eigenspaces with respect to $J,$ that is,
$\cald_\bbc= \cald^{1,0}\oplus \cald^{0,1}.$ Similarly, we get a
splitting of the complexification of the sheaf $\grL^1_B$ of basic
one forms on $M,$ namely
$$\grL^1_B\otimes \bbc = \grL^{1,0}_B\oplus \grL^{0,1}_B.$$
We let $\cale^{p,q}_B$ denote the sheaf of germs of basic forms of
type $(p,q),$ and we obtain a splitting
\begin{equation}
\grL^r_B\otimes \bbc = \bigoplus_{p+q=r}\cale^{p,q}_B.
\end{equation}

The basic cohomology groups $H^{p,q}_B(\calf_\xi)$ are fundamental
invariants of a Sasakian structure which enjoy many of the same
properties as the ordinary Dolbeault cohomology of a K\"ahler
structure. The following theorem was given in \cite{BGN03c}:

\begin{thm} \label{sasbasiccoh} Let $(M,\xi,\eta,\Phi,g)$ be a
compact Sasakian manifold of dimension $2n+1.$ Then we have
\begin{enumerate}
\item $H^{n,n}_B(\calf_\xi)\approx \bbr.$
\item The class $[d\eta]_B\neq 0$ lies in $H^{1,1}_B(\calf_\xi)$ .
\item $H^{p,p}_B(\calf_\xi)>0.$
\item $H^{2p+1}_B(\calf_\xi)$ has even dimension for $p<[n/2].$
\item $H^1(M,\bbr)\approx H^1_B(\calf_\xi).$
\item (Transverse Hodge Decomposition) $H^r_B(\calf_\xi)= \bigoplus_{p+q=r}
H^{p,q}_B(\calf_\xi).$
\item Complex conjugation induces an anti-linear isomorphism
$$H^{p,q}_B(\calf_\xi)\approx H^{q,p}_B(\calf_\xi).$$
\item (Transverse Serre Duality) There is an isomorphism
$$H^{p,q}(\calf_\xi)\approx H^{n-p,n-q}(\calf_\xi).$$
\end{enumerate}
\end{thm}

Theorem \ref{sasbasiccoh} gives rise to fundamental invariants \cite{BGN03b} of
the Sasakian structure, namely the {\it basic Betti numbers} and {\it
basic Hodge numbers}:
\begin{equation}\label{basichodge}
b^B_r(\calf_\xi) =\hbox{dim}~H^{r}_B(\calf_\xi) \qquad
h^{p,q}_B(\calf_\xi)= \hbox{dim}~H^{p,q}_B(\calf_\xi)
\end{equation}
which are related to the basic Betti numbers by
\begin{equation}
b^B_r(\calf_\xi)=\sum_{p+q=r}h^{p,q}(\calf_\xi),
\end{equation}
and satisfy
\begin{equation}\label{cchodge}
h^{p,q}_B(\calf_\xi)=h^{q,p}_B(\calf_\xi),\qquad
h^{p,q}_B(\calf_\xi)=h^{n-p,n-q}_B(\calf_\xi).
\end{equation}

For a description of other transverse holomorphic invariants and relationships between
them, we refer the reader to \cite{BGN03c,BGN03a}. Here we only mention
the {\it basic geometric genus} $p_g(\calf_\xi)=h^{0,n}_B(\calf_\xi),$  and the {\it
basic irregularity} $q=q(\cals)=h^{0,1}_B.$ We remark that $q=
\frac{1}{2}b^B_1(\calf_\xi)=\frac{1}{2}b_1(M)$ is actually a
topological invariant by (5) of Theorem \ref{sasbasiccoh}.
In the case $n=2$, things simplify much more, and all the basic Hodge
numbers are given in terms of the three invariants
$q,p_g(\calf_\xi)$ and $h^{1,1}_B(\calf_\xi).$
Recall now

\begin{defn}\label{quasi-regular}
The characteristic foliation $\calf_\xi$ is said to be {\it
quasi-regular}\index{quasi-regular} if there is a positive integer
$k$ such that each point has a foliated coordinate chart $(U,x)$
such that each leaf of $\calf_\xi$ passes through $U$ at most $k$
times. If $k=1$ then the foliation is called {\it
regular}\index{regular}. If $\calf_\xi$ is not quasi-regular, it
is said to be {\it irregular}. We use the term {\it non-regular} to mean quasi-regular but not
regular.
\end{defn}

In the case of a compact quasi-regular Sasakian manifold $M$ we
know that the space of leaves is a compact Riemannian orbifold
$\calz.$ But since the transverse geometry on $M$ is K\"ahler, the
orbifold must be K\"ahler. But actually in the quasi-regular case
more is true. It follows \cite{BG00a} that $M$ is the total space
of a V-bundle over $\calz,$ and the curvature of the connection
form $\eta$ is precisely the pullback of the K\"ahler form on
$\calz.$ Thus, $\calz$ satisfies an orbifold integrality condition
which we now elaborate. This integrality condition ties Sasakian
geometry on compact manifolds to projective algebraic geometry.
Recall that a compact K\"ahler orbifold $(\calz,\gro)$ is called a
{\it Hodge orbifold} if $[\gro]$ lies in $H^2_{orb}(\calz,\bbz)$.
By a theorem of Baily \cite{Bai57} a Hodge orbifold is a polarized
projective algebraic variety. We have

\begin{thm}\label{quasiregsasastructure}
Let $\cals=(\xi,\eta,\Phi,g)$ be a compact quasi-regular Sasakian structure on a smooth
manifold of dimension $2n+1$, and let $\calz$ denote the space of leaves of the
characteristic foliation. Then
\begin{enumerate}
\item The leaf space $\calz$ is a Hodge orbifold with
K\"ahler metric $h$ and K\"ahler form $\gro$ which defines an
integral class $[\gro]$ in $H^2_{orb}(\calz,\bbz)$ in such a way
that $\pi:(\cals,g) \ra{1.3} (\calz,h)$ is an orbifold Riemannian
submersion. The fibers of $\pi$ are totally geodesic submanifolds
of $\cals$ diffeomorphic to $S^1.$
\item $\calz$ is also a $\bbq$-factorial, polarized, normal projective algebraic variety.
\end{enumerate}
\end{thm}

Conversely, given a Hodge orbifold $(\calz,\gro)$ one can
construct a Sasakian structure on the total space of the circle
bundle defined by the integral class $[\gro]$ in
$H^2_{orb}(\calz,\bbz).$

\bigskip
\section{Deformation Classes of Sasakian Structures}
\medskip
Sasakian structures like K\"ahler structures are plentiful; in fact, the `moduli space' is infinite
dimensional. It is thus important to look for distinguished Sasakian metrics in a given
deformation class, and $\eta$-Einstein metrics provide good examples. There are two
essential types of deformations of Sasakian structures, those that deform the characteristic
foliation (deformations of {\it type I}) and those that fix the characteristic foliation but
deform the contact 1-form (deformations of {\it type II}).

\bigskip
\noindent{\it Deformations of Type I}
\medskip

Let us consider deformations of a Sasakian structure with a fixed transverse
holomorphic structure, that is deformations that deform the characteristic foliation $\calf_\xi$
by deforming the vector field $\xi.$ This type of deformation was first considered by
Takahashi
\cite{Tak} who considered deformations by adding an infinitesimal automorphism $\grr \in
\ga\gu\gt(\xi,\eta,\Phi,g)$. However, here we consider the more general deformations
introduced in \cite{GaOr98}, and used in
the 3-dimensional case by Belgun \cite{Bel01}. Explicitly, we consider deformations of the
contact form whose Reeb vector field is an infinitesimal automorphism of the CR structure,
that is, we consider a positive function $f$ and tensor fields
\begin{equation}\label{Beldef}
\teta=f\eta, \qquad \txi=\xi+\grr,\qquad \tPhi  =\Phi-\Phi\txi\otimes \teta
\end{equation}
such that
\begin{equation}\label{Beldef2}
\teta(\txi)=1,\qquad \txi\rfloor d\teta=0,\qquad \pounds_{\txi}\tPhi=0.
\end{equation}
It is clear from Equations \ref{Beldef} that $\teta$ is a contact 1-form, the contact
subbundle
$\cald=\ker~\eta =\ker~\teta$ remains
unchanged, and that $\tPhi\txi=0.$ Furthermore, if $X$ is a section of $\cald,$ then
$$\tPhi^2X=\tPhi(\Phi X-\teta(X)\Phi\varrho)=\Phi(\Phi X-\teta(\Phi
X)\Phi\varrho)=\Phi^2X=-X.$$
Thus, $\tPhi^2=-\bbi+\txi\otimes \teta,$ so $(\txi,\teta,\tPhi)$ is an underlying almost
contact
structure of the contact structure. Notice that $f$ and $\grr$ are related by
$$f=\frac{1}{1+\eta(\grr)}.$$
Next define the Riemannian metric $\tg$ by
$$\tg=d\teta\circ (\tPhi\otimes \bbi) \oplus \teta\otimes \teta.$$
This metric is compatible with the contact structure and $\txi$ is a Killing vector field of
$\tg.$
Thus, $(\txi,\teta,\tPhi,\tg)$ is K-contact. In fact, we claim that it is Sasakian. To see this
we
need to check normality. But according to Proposition 3.2 of \cite{BG01a} the structure
$(\txi,\teta,\tPhi,\tg)$  is normal if and
only if the almost CR-structure $\tPhi|_\cald$ is integrable and $\tPhi$ is invariant under
$\txi.$ The latter conditions both hold by Equations \ref{Beldef} and \ref{Beldef2}. We have
arrived at

\begin{thm}\label{takdef}
Let $(\xi,\eta,\Phi,g)$ be a Sasakian structure on $M.$ Then the deformations of type I are
deformations through Sasakian structures.
\end{thm}

The metrics $\tg$ and $g$ are related by
$$\tg =\frac{g-\teta \otimes \txi \rfloor g-\txi \rfloor g\otimes \teta
+g(\txi,\txi)\teta\otimes
\teta}{1+ \eta(\rho))} +  \teta\otimes \teta. $$

A very special case of type I deformations occur when $\grr=c\xi$ for some constant $c.$
This does not deform the characteristic foliation, but only rescales the vector field $\xi$ and
contact 1-form $\eta.$ It is convenient to write this transformation in the form
$$\xi'=a^{-1}\xi,\qquad \eta'=a\eta,\qquad \Phi'=\Phi,\qquad
g=ag+a(a-1)\eta\otimes\eta,$$
where $a$ is constant. These were called {\it $\cald$-homothetic deformations} by Tanno
\cite{Tan68}. We have also referred to them as {\it transverse homotheties}.

\bigskip
\noindent{\it Deformations of Type II}
\medskip

Now fix a Sasakian structure $(\xi,\eta,\Phi,g)$ on $M^{2n+1}$ and consider deformations
which stabilize the characteristic foliation
$\calf_\xi.$ More explicitly we wish to
deform $(\xi,\eta,\Phi,g)$ through Sasakian structures that have the same fundamental
basic
cohomology class up to a scale. As in K\"ahler geometry this deformation space is infinite
dimensional. We consider a
deformation of the structure $(\xi,\eta,\Phi,g)$ by adding to $\eta$ a
continuous one parameter family of 1-forms $\grz_t$ that are basic with
respect to the characteristic foliation. We require that the 1-form
$\eta_t=\eta +\grz_t$  satisfy the conditions
\begin{equation}\label{def1}
\eta_0=\eta, \qquad \grz_0=0,\qquad \eta_t\wedge (d\eta_t)^n\neq 0\quad
\forall~~ t\in [0,1].
\end{equation}
This last non-degeneracy condition implies that $\eta_t$ is a
contact form on $M$ for all $t\in [0,1]$ which by Gray's Stability
Theorem \cite{Gra59} belongs to the same underlying contact
structure as $\eta.$ Moreover, since $\grz_t$ is basic $\xi$ is
the Reeb (characteristic) vector field associated to $\eta_t$ for
all $t.$ Now let us define
\begin{align}\label{def2}
\Phi_t&=\Phi -\xi\otimes \grz_t\circ \Phi \\
           g_t&=d\eta_t\circ(\Phi\otimes \bbi_t)+\eta_t\otimes \eta_t.
\end{align}
In \cite{BG03p} it was proved that for all $t\in [0,1]$ and every
basic 1-form $\grz_t$ such that $d\grz_t$ is of type $(1,1)$  and
such that (\ref{def1}) holds $(\xi,\eta_t,\Phi_t,g_t)$ defines a
continuous 1-parameter family of Sasakian structures on $M$
belonging to the same underlying contact structure as $\eta.$

\begin{thm} \label{deftheorem} Let $(M,\xi,\eta,\Phi,g)$ be a
Sasakian manifold. Then for all $t\in [0,1]$ and every basic
1-form $\grz_t$ such that $d\grz_t$ is of type $(1,1)$  and such
that (\ref{def1}) holds $(\xi,\eta_t,\Phi_t,g_t)$ defines a
Sasakian structure on $M$ belonging to the same underlying contact
structure as $\eta.$
\end{thm}

Given a Sasakian structure $\cals =(\xi,\eta,\Phi,g)$ on a manifold
$M,$ we defined $\gF(\xi)$ \cite{BGN03c} to be the family of all Sasakian structures
obtained
by deformations of type II. Notice that any two Sasakian structures $(\xi,\eta,\Phi,g)$ and
$(\xi',\eta',\Phi',g')$ in $\gF(\xi)$ are {\it homologous} in the sense that
$[d\eta']_B=[d\eta]_B.$ This definition was then extended to include transverse
homotheties in \cite{BGN03b}. We defined $\gF(\calf_\xi)$ to be the set
of all Sasakian structures whose characteristic foliation is $\calf_\xi.$
Clearly, we have $\gF(\xi)\subset \gF(\calf_\xi).$ For any Sasakian structure
$\cals=(\xi,\eta,\Phi,g)$ there is the ``conjugate Sasakian structure''
defined by $\cals^c=(\xi^c,\eta^c,\Phi^c,g)=(-\xi,-\eta,-\Phi,g)\in
\gF(\calf_\xi).$ So fixing $\cals$ we define
$$\gF^+(\calf_\xi)=\bigcup_{a\in \bbr^+}\gF(a^{-1}\xi), $$
and $\gF^-(\calf_\xi)$ to be the image of $\gF^+(\calf_\xi)$ under
conjugation. It is then easy to see that we have a decomposition
$\gF(\calf_\xi)=\gF^+(\calf_\xi)\sqcup \gF^-(\calf_\xi).$ Two
Sasakian structures $\cals =(\xi,\eta,\Phi,g)$ and $\cals'
=(\xi',\eta',\Phi',g')$ on a smooth manifold $M$ are said to be
{\it $a$-homologous} if there is an $a\in \bbr^+$ such that
$\xi'=a^{-1}\xi$ and $[d\eta']_B=a[d\eta]_B.$

More generally it is convenient to include both deformations of type I and II as well as
deformations of the CR-structure. For a fixed Sasakian structure $\cals=(\xi,\eta,\Phi,g)$
we denote by $\gF(\cals)$ the space of all Sasakian structures that can be obtained from
$\cals$ by deformations through CR-structures that are compatible with the contact
structure or deformations of type I and II above. The topology of $\gF(\cals)$ is compact
open $C^\infty$-topology induced by the corresponding tensor fields. The underlying
contact structure remains fixed in $\gF(\cals).$

\bigskip
\section{$\eta$-Einstein Metrics and the Calabi Problem}
\medskip
In this section we begin the main objects of study, namely $\eta$-Einstein metrics which
were introduced and studied by Okumura \cite{Oku62} in 1962. In particular, he studied the
relation between the existence of $\eta$-Einstein metrics and certain harmonic forms,
although the correct cohomological setting and connection to the Calabi problem was not
realized until much later.

\begin{defn}\label{etaEinstein}
A Sasakian structure $\cals=(\xi,\eta,\Phi,g)$ on $M$ is said to
be {\it Sasakian $\eta$-Einstein} or just {\it $\eta$-Einstein} if
there are constants $\lambda,\nu$ such that $\Ric_g= \lambda
g+\nu\eta\otimes \eta.$
\end{defn}

One can obviously generalize this definition to any almost contact
metric manifold $(M,\eta,\xi,\Phi,g)$. In such generality one
typically asks for $\lambda,\nu\in C^\infty(M)$ to be arbitrary
functions. However, on a Sasakian or even only K-contact manifold
the Ricci curvature tensor satisfies ${\rm
Ric}_g(\xi,X)=2n\eta(X)$ and it easily follows \cite{Oku62} that

\begin{lem}\label{pos8} Let $M$ be a K-contact manifold of dimension
$2n+1$ such that
$$\Ric_g=\lambda g+\nu\eta\otimes \eta$$
for some $\grl,\nu\in C^\infty(M).$
Then if $n>1$, $\lambda,$ and $\nu$ are constants satisfying $\lambda+\nu=2n$.
\end{lem}

\begin{rem} In the 3-dimensional case any K-contact manifold is automatically Sasakian.
A Sasakian 3-manifold $M$ is always $\eta$-Einstein in the wider
sense of the definition which allows $\lambda,\nu$ to be functions
on $M$. Hence, in this paper  we define a Sasakian 3-manifold to
be $\eta$-Einstein as in Definition \ref{etaEinstein}. On the
other hand, when $M$ is a metric contact manifold one can study
metrics for which $\Ric_g= \lambda g+\nu\eta\otimes \eta$, with
$\lambda,\nu\in C^\infty(M)$ (see \cite{Bl02}).
\end{rem}

An immediate consequence of the definition is

\begin{cor} Every Sasakian $\eta$-Einstein manifold is of
constant scalar curvature $s=2n(\lambda+1)$.
\end{cor}

We briefly recall the basic first Chern class from \cite{BGN03c}.
Now the contact subbundle $\cald$ is a
complex vector bundle and thus has a first Chern class $c_1(\cald)\in
H^2(M,\bbz).$
Consider the long exact sequence \ref{exact} together with the natural map
$H^2(M,\bbz)\ra{1.6} H^2(M,\bbr)$ whose kernel is the torsion part of
$H^2(M,\bbz).$ From (5) of Theorem \ref{sasbasiccoh} we have
\begin{equation}\label{exact2}
\begin{matrix}
          &H^2(M,\bbz)& \\
          &\decdnar{}& \\
          &0\ra{1.2} \bbr\fract{\grd}{\ra{1.5}}
H^2_B(\calf_\xi)\fract{\gri_*}{\ra{1.5}}H^2(M,\bbr)\ra{1.6}
H^1(M,\bbr)\ra{1.6}\cdots.
\end{matrix}
\end{equation}
As in (\ref{exact}) the map $\grd$ is given by $\grd(c)=c[d\eta]$
where $c\in \bbr.$  Now on a Sasakian manifold the vector bundle
$\cald^{1,0}$ is holomorphic with respect to the CR-structure, so
we can compute the free part of $c_1(\cald)=c_1(\cald^{1,0})$ from
the transverse K\"ahler geometry in the usual way. That is
$c_1(\cald)$ can be represented by a basic real closed
$(1,1)$-form $\rho_B.$ The class $c_1^B=[\grr_B]\in
H^2_B(\calf_\xi)$ is independent of the transverse metric and
basic connection used to compute it, and depends only on the
foliated manifold $(M,\calf_\xi)$ with its CR-structure. It is
described in \cite{ElK} and called the {\it basic first Chern
class of $\cald$} there.

\begin{defn}\label{c1def}
A Sasakian structure $(\xi,\eta,\Phi,g)$ is said to be {\it
positive (negative)} if $c_1^B$ is represented by a positive
(negative) definite $(1,1)$-form. If either of these two
conditions is satisfied $(\xi,\eta,\Phi,g)$ is said to be {\it
definite}, and otherwise $(\xi,\eta,\Phi,g)$ is called {\it
indefinite}. $(\xi,\eta,\Phi,g)$ is said to be {\it null} if
$c_1^B=0.$
\end{defn}

In analogy with common terminology of smooth algebraic varieties
we see that a positive Sasakian structure is a {\it transverse
Fano structure}, while a null Sasakian structure is a {\it
transverse Calabi-Yau structure}. The negative Sasakian case
corresponds to the canonical bundle being ample; we refer to this
as a {\it transverse canonical structure}.

Recall the {\it transverse Ricci tensor} $\Ric^T_g$ of $g_T$ to be
the Ricci tensor of the transverse metric $g_T.$ It is related to
the Ricci tensor $\Ric_g$ of $g$
\begin{equation}\label{Riccicurvature}
{\rm Ric}(X,Y)={\rm Ric}_T(X,Y)-2g(X,Y).
\end{equation}
Now as usual define the {\it Ricci form} $\rho_g$
and {\it transverse Ricci form} $\rho^T_g$ by
\begin{equation}\label{trans10b}
\rho_g(X,Y)=\Ric_g(X,\Phi Y), \qquad \rho^T_g(X,Y) =
\Ric^T_g(X,\Phi Y)
\end{equation}
for smooth sections $X,Y$  of $\cald.$ It is easy to check that
these are anti-symmetric of type $(1,1)$ and \ref{Riccicurvature}
implies that they are related by
\begin{equation}\label{trans10c}
\rho^T_g=\rho_g +2d\eta.
\end{equation}
Thus, as in the usual case we have

\begin{lem}\label{c1lemma}
The basic class $2\pi c_1^B\in H^{1,1}_B(\calf_\xi)$ is
represented by the transverse Ricci form $\grr^T_g.$
\end{lem}

In \cite{BGN03c} we used El Kacimi-Alaoui's `Transverse Yau Theorem' \cite{ElK} to prove
the following converse to Lemma \ref{c1lemma}:

\begin{thm} \label{CYSas}
Let $(M,\xi,\eta,\Phi,g)$ be a
Sasakian manifold whose basic first Chern class $c_1^B$ is represented by
the real basic $(1,1)$ form $\grr,$ then there is a unique Sasakian structure
$(\xi,\eta_1,\Phi_1,g_1)\in\gF(\xi)$ homologous to $(\xi,\eta,\Phi,g)$ such
that $\grr_{g_1}=\grr-2d\eta_1$ is the Ricci form of $g_1.$
\end{thm}

In Problem \ref{SCprob} below we give another version of the
Sasakian Calabi Problem. We are interested in certain invariance
properties of $c_1^B$. Since the transverse Ricci form $\grr^T_g$
is invariant under scaling, we see that $c_1^B$ is independent of
the Sasakian structure in $\gF(\calf_\xi).$ Now from the exact
sequence (\ref{exact2}) we have

\begin{prop}\label{c1torsion}
Let $(\xi,\eta,\Phi,g)$ be a Sasakian structure with underlying
contact bundle $\cald.$ Then $c_1(\cald)$ is a torsion class if
and only if there exists a real number $a$ such that
$c_1^B=a[d\eta]_B.$ In particular,
\begin{enumerate}
\item If $c_1(\cald)$ is a torsion class then every Sasakian structure in $\gF(\cals)$ is either
definite or null.
\item If $\gF(\calf_\xi)$ admits a Sasakian $\eta$-Einstein structure, then $c_1(\cald)$ is a
torsion class.
\end{enumerate}
\end{prop}

Hence, a non-torsion $c_1(\cald)$ is the obstruction to the
existence of a Sasakian $\eta$-Einstein metric. In particular, in
the case when $\lambda>-2$  there is a canonical variation to a
\Se metric \cite{Tan79}, so $c_1(\cald)$ is an
obstruction to the existence of a \Se metric on $M$ in the given
$\xi$-deformation class \cite{BG01a}.
Since $c_1(\cald)$ is an invariant of the complex vector
bundle $\cald,$ (2) gives an obstruction to the type of contact structure that can admit a
compatible Sasakian $\eta$-Einstein metric. However, since $c_1(\cald)$ is related to a
topological invariant, namely the second Stiefel-Whitney class $w_2(M),$ there are
topological obstructions to the existence of Sasakian $\eta$-Einstein metrics. From the
natural splitting
$TM=\cald \oplus L_\xi$
where $L_\xi$ is the trivial real line bundle generated by $\xi,$
one gets
$$w_2(M)=w_2(TM)=w_2(\cald)$$
which is the mod 2 reduction of $c_1(\cald)\in H^2(M,\bbz).$ Since
a manifold $M$ admits a spin structure if and only if $w_2(M)=0,$
Proposition \ref{c1torsion} implies

\begin{thm}\label{spinobs}
Let $M$ be a non-spin manifold (i.e., $w_2(M)\neq 0$) with
$H_1(M,\bbz)$ torsion free. Then $M$ does not admit a Sasakian
$\eta$-Einstein structure.
\end{thm}

Thus, Theorem \ref{spinobs} includes
as a special case the well-known fact \cite{BG00a, Mor97} that a
simply connected \Se manifold is necessarily spin.
Next we give an example of a Sasakian manifold that does not admit a Sasakian
$\eta$-Einstein metric.

\begin{exmp}
Consider the non-trivial $S^3$-bundle over $S^2,$ denoted by $X_\infty$ in Barden
\cite{Bar65}. We represent $X_\infty$ as a non-trivial circle bundle over a Hirzebruch
surface. Recall the Hirzebruch surfaces $S_n$ \cite{GrHa78} are realized as the
projectivizations of the sum of two line bundles over $\bbc\bbp^1,$ namely
$$S_n=\bbp\bigl(\calo \oplus\calo(-n)\bigr).$$
They are diffeomorphic to $\bbc\bbp^1\times \bbc\bbp^1$ if $n$ is
even, and to the blow-up of $\bbc\bbp^2$ at one point, which we
denote as $\widetilde{\bbc\bbp}^2,$ if $n$ is odd. We need to take
$n$ odd to get a non-spin circle bundle. Now ${\rm
Pic}(S_n)\approx \bbz\oplus \bbz,$ and we can take the Poincar\'e
duals of a section of $\calo(-n)$ and the homology class of the
fiber as its generators. The corresponding divisors can be
represented by rational curves which we denote by $C$ and $F,$
respectively satisfying
$$C\cdot C=-n, \qquad F\cdot F=0,\qquad C\cdot F=1.$$
Let $\gra$ and $\grb$ denote the Poincar\'e duals of $C$ and $F$, respectively, and write
the K\"ahler form $\gro$ so that $[\gro]=a\gra +b\grb.$ According to \cite{HwSi97} every
K\"ahler class can be represented in this form with $b=2$ and $a>0.$
Let
$$S^1\ra{1.7}M\fract{\pi}{\ra{1.7}}S_n$$
be the circle bundle over $S_n$ whose first Chern class is $[\gro].$ Since we want $M$
simply connected, we need to choose $\gcd(a,b)=1.$ It suffices to take $n=3,$ and
$[\gro]=\gra+2\grb.$ Now $c_1(\cald)=\pi^*c_1(S_n)=-\pi^*c_1(K)$ where $K$ is the
canonical divisor of $S_n.$ From \cite{GrHa78} we have $K=-2C-(n+2)F.$ So in our case
we have $c_1(K)=-2\gra-5\grb.$ From the Gysin sequence of the circle bundle we have
$$0\ra{1.5}\bbz\fract{\cup~[\gro]}{\ra{2.5}}H^2(S_n,\bbz)\fract{\pi^*}{\ra{2.5}}H^2(M,\bbz)
\ra{1.5}0,$$
giving $c_1(\cald)=\pi^*(c_1(-K))=\pi^*(2\gra+5\grb)=\pi^*\grb.$ Hence, $w_2(M)\neq 0$
implying that $M=X_\infty$ by Barden's Classification Theorem \cite{Bar65}. Thus, we have
produced a Sasakian structure on $X_\infty$, but by Theorem \ref{spinobs} it cannot admit
a Sasakian $\eta$-Einstein structure.
\end{exmp}

We are interested in proving the existence of Sasakian $\eta$-Einstein metrics in a given
deformation class of Sasakian structures.

\begin{prob}\label{SCprob}
{\bf [Sasakian Calabi Problem]} Given a manifold $M$ with Sasakian
structure $\cals=(\xi,\eta,\Phi,g)$ and with basic first Chern
class $c_1^B$ that is represented by either a positive definite,
negative definite real basic $(1,1)$ form $\grr^T,$ or if $c_1^B$
vanishes, does there exist a Sasakian structure $\cals'\in
\gF(\calf_\xi)$ with an $\eta$-Einstein metric $g'$?
\end{prob}

The solution to this problem is, of course, given in local
CR-coordinates $(z_i,\bar{z}_i,x)$ on $M$ by the ``transverse
Monge-Amp\`ere equation''
\begin{equation} \label{mon-amp}
\frac{\det(g^T_{i\barj}+\phi_{i\barj})}{\det(g^T_{i\barj})}=
e^{-k\phi+F}, \quad g^T_{i\barj}+\phi_{i\barj}>0,
\end{equation}
where $g^T$ is the transverse metric, $\phi$ and $F$ are real
basic functions, and $\phi_{i\barj}$ are the components of
$\partial\bar{\partial}\phi =d\grz_t$ with respect to the
transverse coordinates $(z_i,\bar{z}_{\barj}).$  The cases when
$c_1^B$ is zero, positive or negative definite correspond to
finding solutions of equation \ref{mon-amp} with the constant
$k=0,>0,<0,$ respectively. In the negative and zero cases we have
the transverse version of theorems of Yau \cite{Yau78} and Aubin
\cite{Aub82} give

\begin{thm} \label{aub-yau} If the class $[\grr^T]\in
H^2_B(\calf_\xi)$ is zero or can be represented by a negative
definite $(1,1)$ form, then there exists a Sasakian structure
$\cals\in \gF(\xi)$ with an $\eta$-Einstein metric $g$ on $M$
with $\lambda=-2$ in the first case and $\lambda< -2$ in the
second.
\end{thm}

In other words, in complete analogy with the K\"ahlerian case, we see that there are no
obstructions to solving Sasakian or transverse Monge-Amp\`ere
problem in the negative or null case.

\bigskip
\section{Further Results on $\eta$-Einstein Metrics}
\medskip
Let us now consider a Sasakian  $\eta$-Einstein with Sasakian
structure $\cals=(\xi,\eta,\Phi,g)$  and $\eta$-Einstein constants
$(\lambda,\nu)$. A simple calculation shows that
$\cald$-homotheties preserve $\eta$-Einstein condition. More
precisely, we have
\begin{prop}\label{etahomothetic}
Let  $(M,\xi,\eta,\Phi,g)$ be a Sasakian $\eta$-Einstein manifold
with constants $(\lambda,\nu)$. Consider a $\cald$-homothetic
Sasakian structure $\cals'=(\xi',\eta',\Phi',g')=(
a^{-1}\xi,a\eta,\Phi,ag+a(a-1)\eta\otimes\eta)$. Then $(M,\cals')$
is also $\eta$-Einstein with constants
$$\lambda'=\frac{\lambda+2-2a}{a},\qquad \nu'=2n-\frac{\lambda+2-2a}{a}.$$
\end{prop}
\begin{proof}
The proof follows from a simple formula relating Ricci curvature
tensors of the Sasakian metric $g',g$ related by a
$\cald$-homothetic transformation \cite{Tan68}:
$${\rm Ric}_{g'}={\rm Ric}_{g}-2(a-1)g+(a-1)(2n+2+2na)\eta\otimes\eta.$$
\end{proof}
Note that $\cald$-homotheties must preserve the sign of the
Sasakian $\eta$-Einstein structure (invariance of basic Chern
classes) and they are completely analogous to the usual
homotheties of Einstein metrics. When $(M,\xi,\eta,\Phi,g)$ is
null then always $(\lambda,\nu)=(-2,2n+2)$. When
$(M,\xi,\eta,\Phi,g)$ is positive then $\lambda>-2$. But by
applying a suitable $\cald$-homothety one can get any value of
$\lambda\in(-2,\infty)$. Similarly, when $(M,\xi,\eta,\Phi,g)$ is
negative $\lambda<-2$ and by applying a suitable $\cald$-homothety
one can get any value of $\lambda\in(-\infty,-2)$.

An easy computation shows that, in addition to the structure
Theorem \ref{quasiregsasastructure}, in the $\eta$-Einstein case
we can actually say more.

\begin{thm}\label{etasasastructure}
Let $(M,\xi,\eta,\Phi,g)$ be a compact quasi-regular Sasakian
$\eta$-Einstein manifold of dimension $2n+1$, and let $\calz$
denote the space of leaves of the characteristic foliation. Then
the leaf space $\calz$ is a Hodge orbifold with K\"ahler-Einstein
metric $h$ with Einstein constant $\lambda+2$.
\end{thm}

The case of $\lambda+2>0$ is special as follows from Proposition
\ref{etahomothetic}. By applying suitable $\cald$-homothety one
can always choose $\lambda=2n$ which turns the $\eta$-Einstein
metric into an Einstein metric.  This observation is originally
due to Tanno who used it to show that a unit tangent bundle of
$S^n$ has a homogeneous Sasakian-Einstein structure. When $n=3$
one gets a homogeneous Sasakian-Einstein metric on $S^2\times S^3$
\cite{Tan79}. We have

\begin{thm}\label{se}
Let $(M,\xi,\eta,\Phi,g)$ be a compact quasi-regular Sasakian
$\eta$-Einstein manifold of dimension $2n+1$, and let
$\lambda>-2$. Let $(\calz,h)$ denote the space of leaves of the
characteristic foliation with its K\"ahler-Einstein metric $h$.
Then $M$ admits a Sasakian-Einstein structure
$(\xi',\eta',\Phi',g')\in \gF(\calf_\xi)$, $\cald$-homothetic to
$\cals$, with  $g'=\alpha g+\alpha(\alpha-1)\eta\otimes\eta$ and
$\alpha =\frac{\lambda+2}{2n+2}.$
\end{thm}

It is useful therefore to think of positive Sasakian
$\eta$-Einstein structure as a family parameterized by a positive
constant $a\in\bbr^+$. In addition, if we choose the
Sasakian-Einstein structure in this family $\cals_1$ as a
reference point we will distinguished between two subfamilies.

\begin{defn}\label{squashed}
Let $(M,\xi_1,\eta_1,\Phi_1,g_1)$ be a compact quasi-regular
Sasakian-Einstein manifold of dimension $2n+1$, and let
$\cals_a=(\xi_a,\eta_a,\Phi_a,g_a)$ be the Sasakian
$\eta$-Einstein structure obtained from $\cals_1$ by applying
$\cald$-homothetic transformation to $\cals_1$ with constant
$a\in\bbr^+$. We say that $\cals_a$ is {\it stretched} when $a>1$
$(-2<\lambda<2n, \nu>0)$ and {\it squashed} when $0<a<1$
$(\lambda>2n, \nu<0 )$.
\end{defn}

Note that positive squashed Sasakian $\eta$-Einstein manifolds can
have scalar curvature of any sign. There is a unique metric in the
family $\cals_a$ which makes $\lambda=-1$ in which case the
Sasakian manifold becomes scalar-flat.

In the negative case there is another interesting structure
related to Sasakian $\eta$-Einstein metric. Recall \cite{Bau00, Boh03}
the following

\begin{defn}\label{squashed2}
Let $(M,g)$ be a Lorentzian manifold of dimension $2n+1$ and let
$\xi$ be a time-like Killing vector field such that
$g(\xi,\xi)=-1$. We say that that $M$ is Sasakian if
$\Phi(X)=\nabla_X\xi$ satisfies the condition
$(\nabla_X\Phi)(Y)=g(X,\xi)Y+g(X,Y)\xi$ and Sasakian-Einstein if,
in addition, $g$ is Einstein.
\end{defn}
Equivalently, we can require the metric cone $\calc(M)=
(\bbr^+\times M,-dt^2+r^2g, d(r^2\eta))$ to be pseudo-K\"ahler of
signature $(2,2n)$. In the Sasakian-Einstein case the cone metric
is pseudo-Calabi-Yau and the Einstein constant must equal to
$-2n$. This is in complete analogy with the Riemannian case.
Geometries of this type are interesting as they provide examples
of the so-called twistor spinors on Lorentzian manifolds. As a
simple consequence we get

\begin{prop}\label{lorse}
Let $(M,\xi,\eta,\Phi,g)$ be a compact quasi-regular Sasakian
$\eta$-Einstein manifold of dimension $2n+1$, and let
$\lambda<-2$. Let $(\calz,h)$ denote the space of leaves of the
characteristic foliation with its K\"ahler-Einstein metric $h$.
Then $M$ admits a Lorentzian Sasakian-Einstein structure such that
$$\xi'=a^{-1}\xi,\qquad -g'=ag+a(a-1)\eta\otimes\eta,$$
where $a=\frac{\lambda+2}{2+2n}.$\end{prop}

Combining this with Theorem \ref{aub-yau} we have

\begin{cor} Every negative Sasakian manifold
admits a Lorentzian Sasakian-Einstein structure.
\end{cor}

\bigskip
\section{$\eta$-Einstein Structures on Links}
\medskip
Let $\bfw=(w_0,\ldots,w_n)$ be any positive vector in
$\bbr^{n+1},$ that is $w_i>0$ for all $i=0,\ldots,n.$ The vector
field corresponding to $\bfw$ is denoted by $\xi_\bfw$, and the
corresponding 1-form by $\eta_\bfw.$ Let
$(x_0,\ldots,x_n,y_0,\ldots,y_n)$ be the coordinates on the unit
sphere $S^{2n+1}$ in $\bbr^{2n+2}$. A family Sasakian structures
$\cals_\bfw=(\xi_\bfw,\eta_\bfw,\Phi_\bfw,g_\bfw)$ is now defined
by
\begin{equation}\label{weightedxieta}
\xi_\bfw
=\sum_{i=0}^nw_i(y_i\partial_{x_i}-x_i\partial_{y_i}),\qquad
\eta_\bfw =\frac{\sum_{i=0}^nx_idy_i}{\sum_{i=0}^nw_i((x^i)^2
+(y^i)^2)},
\end{equation}
with $\Phi_\bfw$ and $g_\bfw$ further determined by $\xi_\bfw$ and
$\eta_\bfw$. We shall refer to
$(\xi_\bfw,\eta_\bfw,\Phi_\bfw,g_\bfw)$ as the {\it weighted
Sasakian structure} on $S^{2n+1},$ and we denote the unit sphere
with this Sasakian structure by $S^{2n+1}_\bfw.$ Note that the
standard Sasakian-Einstein structure (Hopf fibration) corresponds
to setting $\bfw=(1,\ldots,1)$. The deformed structures
$(\xi_\bfw,\eta_\bfw,\Phi_\bfw,g_\bfw)$ are not all distinct. The
Weyl group $\calw(SU(n+1))={\rm Nor}(\gT_{n+1})/\gT_{n+1}$ acts as
outer automorphisms on the maximal torus, and it is well-known
that $\calw(SU(n+1))=\grS_{n+1}$, the permutation group on $n+1$
letters. So for any $\grs\in \grS_{n+1}$ the Sasakian structures
$\cals_\bfw=(\xi_\bfw,\eta_\bfw,\Phi_\bfw,g_\bfw)$ and
$(\xi_{\grs(\bfw)},\eta_{\grs(\bfw)},\Phi_{\grs(\bfw)},g_{\grs(\bfw)})$
are isomorphic. Thus, we shall frequently take the weights to be
ordered such that $w_0\leq w_1\leq \cdots \leq w_n.$ It is
important to note that all the Sasakian structures $\cals_\bfw$
belong to the same underlying contact structure on $S^{2n+1},$
namely the standard one.

Consider the affine space $\bbc^{n+1}$ together with a weighted
$\bbc^*$-action given by
\begin{equation}\label{C*action}
(z_0,\ldots,z_n)\mapsto
\grl\cdot\bfz=(\grl^{w_0}z_0,\ldots,\grl^{w_n}z_n),
\end{equation}
where $\grl\in \bbc^*,$ and the {\it weights} $w_j$ are positive
integers.

\begin{defn}\label{whp}
A polynomial $f\in \bbc[z_0,\ldots,z_n]$ is said to be {\it
weighted homogeneous} of {\it degree} $d$ if there are positive
integers $w_0,w_1,\ldots,w_n,d$ such that
$$f(\grl\cdot\bfz)=f(\grl^{w_0}z_0,\ldots,\grl^{w_n}z_n)=\grl^df(z_0,\ldots,z_n)=\grl^df(\bfz).$$
\end{defn}

We shall assume that $w_0\leq w_1\leq \cdots\leq w_n$ and that
$\gcd(w_0,\ldots,w_n)=1$ unless otherwise stated. The zero set
$V_f$ of a non-degenerate weighted homogeneous polynomial $f,$ is
just the {\it weighted affine cone} $C_f.$ If the origin is an
isolated singularity in $C_f,$ then it is the only singularity of
$C_f.$ We shall assume this from now on. In this case we can take
$\gre=1$ so that $S^{2n+1}_\gre(\bfz_0)$ becomes the unit sphere
$S^{2n+1}$ centered at the origin with its weighted Sasakian
structure $\cals_\bfw$. In this case our link
\begin{equation}\label{isolink}
L_f= C_f\cap S^{2n+1}.
\end{equation}
is a smooth manifold of dimension $2n-1$ which by the Milnor
Fibration Theorem is $(n-2)$-connected.

One can easily see that every link is either positive, negative, or
null. More precisely, we have \cite{BGK03}
\begin{thm}\label{saslink}
Let $f$ be a non-degenerate weighted homogeneous polynomial of
degree $d$ and weight vector $\bfw.$ The Sasakian structure
$(\xi_\bfw,\eta_\bfw,\Phi_\bfw,g_\bfw)$ on $S^{2n+1}$ induces by
restriction a Sasakian structure, denoted by
$\cals_{\bfw,f}=\cals_\bfw\mid_{L_f}$ on the link $L_f.$
Furthermore, let $|\bfw|=w_0+\cdots+w_n$. Then
\begin{enumerate}
\item  $\cals_{\bfw,f}$ is positive when $d-|\bfw|<0$,
\item  $\cals_{\bfw,f}$ is null when $d-|\bfw|=0$,
\item  $\cals_{\bfw,f}$ is negative when $d-|\bfw|>0$.
\end{enumerate}
\end{thm}

\begin{exmp} Let $f(\bfz)=z_0^{a_0}+\cdots+z_n^{a_n}$. Then the link
$L_f$ is called the Brieskorn-Pham link and we will denote it by
$L(\bfa)=L(a_0,\ldots,a_n)$. It can be seen that the weighted
degree of the Brieskorn-Pham polynomial is  $d={\rm
lcm}(a_0,\ldots,a_n)$ and  the weights are $w_j=d/a_j$. Then by
Theorem \ref{saslink}
\begin{enumerate}
\item  $L(\bfa)$ is positive when $\sum_i\frac{1}{a_i}>1$,
\item  $L(\bfa)$ is null when $\sum_i\frac{1}{a_i}=1$,
\item  $L(\bfa)$ is negative when $\sum_i\frac{1}{a_i}<1$.
\end{enumerate}
\end{exmp}

Theorem \ref{saslink} together with Theorem \ref{aub-yau} provide a rich source of
$\eta$-Einstein metrics in the non-positive
case where there are no obstructions. Summarizing we have

\begin{thm}\label{nullnegCal}
Let $\cals_{\bfw,f}$ be null or negative. Then $L_f$ admits a Sasakian
structure $\cals'$ with $\eta$-Einstein metric $g'$ in the same deformation class as
 $\cals_{\bfw,f}.$
\end{thm}

In the Fano case there are obstructions to the existence of such
structures. However, as every $\eta$-Sasakian metric can be
deformed to a Sasakian-Einstein one such obstruction are
completely equivalent to the obstructions preventing existence of
orbifold K\"ahler-Einstein metric on the transverse space
$\calz_f$. Existence of K\"ahler-Einstein metrics on various Fano
orbifolds can be often established by the continuity method (cf.
\cite{BGK03} and references therein).

We finish this section by reviewing very briefly how one determines the topology of the links.
The procedure is that of Milnor and Orlik \cite{MiOr70} and involves computing the
Alexander polynomial of the link. For more details we refer to \cite{MiOr70} and our previous
work \cite{BG01b,BGN03c}. Milnor and Orlik give a combinatorial formula for the
$(n-1)$-st Betti number of a $(2n-1)$-dimensional link, viz.
\begin{equation}\label{Bettiequation}
b_{n-1}(L_f)= \sum (-1)^{n+1-s}\frac{u_{i_1}\cdots u_{i_s}}{v_{i_1}\cdots v_{i_s}{\rm
lcm}(u_{i_1},\ldots,u_{i_s})},
\end{equation}
where the quotients $d/w_i$ are written in the irreducible form $u_i/v_i$, and the sum is
taken over all the $2^{n+1}$ subsets $\{i_1,\ldots,i_s\}$ of $\{0,\ldots,n\}.$
In the case of
determining homology spheres or more generally rational homology spheres it is more
efficient to use the Brieskorn Graph Theorem. We refer to \cite{BGK03} for details. Another
important result proven in \cite{BG01b} concerns 5-dimensional manifolds only. It says that
if the weights $(w_0,w_1,w_2,w_3)$ satisfy the ``well-formedness'' condition
$\gcd(w_i,w_j,w_k)=1$ for all distinct $i,j,k,$ then the 5-dimensional link $L_f$ has no
torsion.

\section{Sasakian Structures in Dimension 3}

Three dimensional Sasakian geometry is well understood, culminating in the recent
uniformization theorem due to Belgun \cite{Bel00} which we state precisely below. Geiges
\cite{Gei97} showed that a compact 3-manifold admits a Sasakian structure if and only if it
is diffeomorphic to one of the following:
\begin{enumerate}
\item $S^3/\grG$ with $\grG\subset \gI_0(S^3)=SO(4).$
\item $\widetilde{SL}(2,\bbr)/\grG$ where is universal cover of $SL(2,\bbr)$ and
$\grG\subset \gI_0(\widetilde{SL}(2,\bbr)).$
\item ${\rm Nil}^3/\grG$ with $\grG\subset \gI_0({\rm Nil}^3).$
\end{enumerate}
Here $\grG$ is a discrete subgroup of the connected component $\gI_0$ of the
corresponding
isometry group with respect to a `natural metric', and ${\rm Nil}^3$ denotes the 3 by 3
nilpotent real matrices, otherwise known as the Heisenberg group. These are three of the
eight model geometries of Thurston \cite{Thu97}, and correspond precisely to the compact
Seifert bundles with non-zero Euler characteristic \cite{Sco83}. For further discussion of the
isometry groups we refer to \cite{Sco83}. We refer to the three model geometries above as
{\it spherical}, $\widetilde{SL_2}$ type, and {\it nil} geometry, respectively.

The purpose of this section then is to describe some examples represented as links of
isolated hypersurface singularities. First note that in dimension 3 the basic cohomology
group
$H^2(\calf_\xi)$ is 1-dimensional so we cannot get indefinite forms
other than zero. Hence

\begin{prop}\label{sas3}
Let $M$ be a compact 3-dimensional Sasakian manifold. Then $M$ is
is either positive, negative, or null.
\end{prop}
These 3 types of Sasakian structures correspond precisely to the 3  model geometries
above.
Belgun's Theorem uniformizes these three cases. We prefer to rephrase Belgun's theorem
in terms of the earlier work of Tanno \cite{Tan69b} on constant $\Phi$-sectional curvature.
In dimension three the $\Phi$-sectional curvature is determined by one function, namely
$H(X)=K(X,\Phi X).$ Thus, in analogy with 3-dimensional Einstein geometry we have (see also
\cite{Gui02})

\begin{prop}\label{3etaEin}
A three dimensional Sasakian manifold is $\eta$-Einstein if and only if it has constant
$\Phi$-sectional curvature.
\end{prop}

\begin{proof}
The if part, which holds in all dimensions, follows easily from Tanno's classification theorem
\cite{Tan69b}. To prove the only if part we choose a local orthonormal bases $X,\Phi X,\xi$
for the Sasakian structure $\cals=(\xi,\eta,\Phi,g).$  Now in three dimensions the Ricci
curvature and the sectional curvature are related by \cite{Pet98}
\begin{align*}
\Ric_g(X,X) &=K(X,\Phi X) +K(X,\xi) =K(X,\Phi X) +1 \\
\Ric_g(\Phi X,\Phi X) &=K(X,\Phi X) +K(\Phi X,\xi) =K(X,\Phi X) +1 \\
\Ric_g(\xi,\xi)&=K(X,\xi)+K(\Phi X,\xi) = 2
\end{align*}
So if $g$ is $\eta$-Einstein we have
$$\Ric_g =\gra g+ (2-\gra)\eta\otimes \eta,$$
and this gives $K(X,\Phi X)=\gra -1.$
\end{proof}

Now if a Sasakian structure $\cals=(\xi,\eta,\Phi,g)$ has constant $\Phi$-sectional
curvature $c$, then it can be transformed to one of the cases $c=1,-3,$ or $-4$ by a
transverse homothety. The  Uniformization Theorem can now be stated as:

\begin{thm}[Uniformization \cite{Bel00}]\label{Belthm}
Let $M$ be a 3-dimensional compact manifold admitting a Sasakian structure
$\cals=(\xi,\eta,\Phi,g).$ Then
\begin{enumerate}
\item If $\cals$ is positive, $M$ is spherical, and there is a Sasakian metric of constant
$\Phi$-sectional curvature $1$ in the same deformation class as $g.$
\item If $\cals$ is negative, $M$ is of $\widetilde{SL_2}$ type, and there is a Sasakian
metric of constant $\Phi$-sectional curvature $-4$ in the same deformation class as $g.$
\item If $\cals$ is null, $M$ is nil, and there is a Sasakian metric of constant
$\Phi$-sectional curvature $-3$ in the same deformation class as $g.$
\end{enumerate}
\end{thm}

In the positive case the universal cover $\tM$ is $S^3$ with its standard round sphere
metric. In this case one generally needs both type I and II deformations. In the negative
case $\tM=B^2\times \bbr$, the product of the unit 2-ball with $\bbr$; whereas, in the null
case $\tM=\bbr^3.$ In both of these cases the metrics have constant $\Phi$-sectional
curvature, but not constant sectional curvature. Moreover, deformations of type I do not
exist in the both negative and null cases except for $\cald$-homotheties. Clearly, these
`standard metrics' are all $\eta$-Einstein. More generally in any dimension constant
$\Phi$-sectional curvature metrics are $\eta$-Einstein.

In the positive case deformations of type I are definitely needed. Theorem \ref{Belthm}
does not hold without them.  A counterexample
is the weighted Sasakian structure on $S^3$ with weights
$\bfw=(p,q)$, $p\not=1$. The transverse space of the
characteristic foliation is the weighted projective line
$\bbp_\bbc^1(p,q)$ which is an example of a Fano orbifold which
does not admit an orbifold  metric of constant curvature. More
generally, the existence of K\"ahler-Einstein metrics on weighted
projective spaces of arbitrary dimension is obstructed by the
non-vanishing of the Futaki character \cite{ACGT04}.

\begin{rem} A Sasakian-Einstein
structure on a 3-Sasakian manifold does not have to be a part of
the 3-Sasakian structure. The simplest example when this is the
case is the lens space $\bbz_k\backslash S^3$. Consider the unit
3-sphere $S^3\simeq Sp(1)$ as the unit quaternion $\sigma\in\bbh$.
Such a sphere has two 3-Sasakian structures generated by the left
and the right multiplication. Consider the homogeneous space
$\bbz_k\backslash S^3$, where the $\bbz_k$-action is given by the
multiplication from the left by $\rho\in Sp(1)$, $\rho^k=1$. The
quotient still has the ``right" 3-Sasakian structure. But it also
has a  ``left" Sasakian structure (the centralizer of $\bbz_k$ in
$Sp(1)$ is an $S^1$ and it acts on the coset from the left). This
left Sasakian structure is actually regular while none of the
Sasakian structures of the right 3-Sasakian structure can be
regular unless $k=1,2$.
\end{rem}

\begin{exmp} Consider the standard metric on $S^3\ra{1.2}S^2$ of constant
sectional curvature equal to 1. This metric is Sasakian-Einstein
and 3-Sasakian. Applying $\cald$-homothety to this metric produces
a family of $\eta$-Einstein U(2)-invariant metrics which can be
written as $\sigma_1^2+\sigma_2^2+a\sigma_3^2$, where $\sigma_1$,
$\sigma_2$ and $\sigma_3$ are the standard left-invariant
$1$-forms on $S^3\simeq SU(2)$ and $a$ is a nonzero constant. The
metrics were first considered by Berger. We will return to this
example in the last section as it it the most basic Einstein-Weyl
manifold.
\end{exmp}

Let us now consider some examples of Sasakian
structures on 3-manifolds represented as links of isolated hypersurface singularities defined
by weighted homogeneous polynomials.

\begin{exmp}\label{posex}Positive Sasakian 3-manifolds.
Belgun's Theorem \ref{Belthm} says that every positive Sasakian 3-manifold is a spherical
space form $S^3/\grG$ where $\grG$ is a finite subgroup of $SU(2).$ These are precisely
\begin{enumerate}
\item $\grG= \bbz_p$ the cyclic group of order $p,$
\item $\grG= \bbd_m^{*}$ a binary dihedral group where $m$ is an integer greater than 2,
\item $\grG= \bbt^{*}$ the binary tetrahedral group,
\item $\grG= \bbo^{*}$ the binary octahedral group,
\item $\grG= \bbi^{*}$ the binary icosahedral group.
\end{enumerate}
Orlik \cite{Or70} shows that all of these groups can be realized by links of weighted
homogeneous polynomials,
and that this exhausts all cases with $|\bfw|>d.$ It should be noted, however, that not all
lens spaces $L(p,q)$ occur, but only $L(p,1).$ The table below indicates the subgroup
$\grG\subset SU(2)$ together with a representative polynomial and the corresponding
Dynkin diagram of the singularity.
\begin{center}\vbox{\[\begin{array}{|l|l|l|} \hline {\rm Diagram}&{\rm Subgroup}
&{\rm Polynomial}\hl\hline
A_{p-1}&\bbz_p&z_0^p+z_1^2+z_2^2\hl
D_m&\bbd_m^*&z_0^2z_1+z_1^m+z_2^2\hl
E_6&\bbt^*&z_0^4+z_1^3+z_2^2\hl
E_7&\bbo^*&z_0^3+z_1^3z_0+z_2^2\hl
E_8&\bbi^*&z_0^5+z_1^3+z_2^2\hl
\end{array}\]}
\end{center}
Notice that $L(5,3,2)=S^3/\bbi^*$ is the famous Poincar\'e homology sphere.
\end{exmp}

\begin{exmp}\label{null3man}Null Sasakian 3-manifolds.
Up to covering these are all non-trivial circle bundles over a complex torus, and it is known
\cite{Mil75} that all such bundles can be constructed from the three dimensional Heisenberg
group $H\simeq \bbr^3$ which can be represented as a subgroup of $GL(3,\bbr)$ by
considering matrices of the form
\begin{equation} \bba=\left(
\begin{array}{ccc}
1 & x & z\\
0 & 1 & y\\
0 & 0 & 1 \end{array}\right)\  ;\ x,y,z\in\bbr.
\end{equation}
We also define subgroups $H_k\subset H$ of matrices for which
$x,y,z$ are integers divisible by some positive integer $k$. It is
easy to see that the quotient manifold $M_k=H/H_k$ is a circle
bundle over a torus with Chern number $\pm k$. The fibration $\pi
:M_k\ra{1.2} T^2 $ is explicitly given by
\begin{equation}
\pi(\bba) =(x\ {\rm mod}\ k, z\ {\rm mod}\ k).
\end{equation}
Now Sasaki \cite{Sas65} showed that $\bbr^3$ has a natural null Sasakian $\eta$-Einstein
metric with constants $(\lambda,\nu)=(-2,4),$
although the connection with the Heisenberg group wasn't observed until later. This metric is
\begin{equation}\label{metr}
g=dx^2+dy^2 +(dz-y~dx)^2.
\end{equation}
The Sasakian structure
$\cals=(\xi,\eta,\Phi, g)$ on $H$ is defined by the formulas
\begin{equation} \left\{
\begin{array}{lcl}
\eta=dz-y\ dx & \xi =\frac{\partial}{\partial z} & \\
\Phi =(\frac{\partial}{\partial x}+y\frac{\partial}{\partial z})\otimes dy
-\frac{\partial}{\partial y}\otimes dx.
\end{array}\right.
\end{equation}
Now one can easily check that the tensor fields defining $\cals$ are invariant under the
action of $H$ on itself. Hence, $\cals$ defines a Sasakian structure on the compact
quotient manifolds $M_k$.
Milnor \cite{Mil75} observed that the integral homology can easily be computed, viz.
$$H_1(M_k,\bbz)=H_k/[H_k,H_k]=\bbz\oplus\bbz\oplus\bbz_k.$$
However, as indicated by the table below only the cases $k=1,2,3$ can be realized as links
of isolated hypersurface singularities of weighted homogeneous polynomials.
There are precisely three null Sasakian 3-manifolds that can be represented by links of
isolated hypersurface singularities of weighted homogeneous polynomials \cite{Or70}.
These are realized by arbitrary weighted homogeneous polynomials of degrees $6,4$ and
$3$, and they exhaust all possibilities that satisfy $|\bfw|=d$. The table below gives the
weights, degrees, and number of monomials of the three polynomials.
\begin{center}\vbox{\[\begin{array}{|c|c|c|c|} \hline {\rm manifold}& {\rm weight ~vector}
& {\rm degree} &{\rm \#~of~monomials}\hl\hline
M_1 &(1,2,3)&6&7\hl
M_2 &(1,1,2)&4&9\hl
M_3 &(1,1,1)&3&10\hl
\end{array}\]}
\end{center}
\end{exmp}

\begin{exmp}\label{negsas}Negative Sasakian 3-manifolds.
These are up to covering circle bundles over Riemann surfaces of genus $g>1.$ There are
infinitely many examples which can be represented as links of isolated hypersurface
singularities defined by weighted homogeneous polynomials. Milnor
proves \cite{Mil75} that in the case of Brieskorn 3-manifolds
$L(a_0,a_1,a_2)=\widetilde{SL}(2,\bbr)/\Gamma$, where
$\widetilde{SL}(2,\bbr)$ is the universal cover of $SL(2,\bbr)$ and
$\Gamma$ is a discrete subgroup completely determined by the
triple $(a_1,a_2,a_3)$. Perhaps the most interesting are the integral homology spheres. In
this case a theorem of Saeki \cite{Sae87} says that the link of an isolated hypersurface
singularity
defined by a weighted homogeneous polynomial has the same knot type as a Brieskorn link.
Now Breiskorn's graph Theorem \cite{Bri66} implies that a Brieskorn link $L(\bfa)$ is a
homology sphere if and only if
$\gcd(a_i,a_j)=1$ for all $i\neq j=0,1,2.$ So clearly 3-dimensional homology
spheres are quite numerous, and all but one, the Poincar\'e homology sphere $L(5,3,2),$
have
infinite fundamental group. More, generally for any link of an isolated hypersurface singularity
defined by a weighted homogeneous polynomial, Orlik \cite{Or70} proved that $\pi_1$ is
infinite if and only if $d\geq |\bfw|,$ and $\pi_1$ is infinite nilpotent if and only if $d=|\bfw|$
which is the Euclidean or null Sasakian case. Obviously, for integral homology spheres
$\pi_1$ must be perfect, i.e. $[\pi_1,\pi_1]=\pi_1$ where $[\pi_1,\pi_1]$ denotes the
commutator subgroup. A rather large class of examples can be treated by considering
Brieskorn polynomials of the form
$$f=z_0^{6k\pm 1}+z_1^3+z_2^2.$$
Since $d-|\bfw|=6k-6\pm 1\geq
0,$ $\pi_1(L(6k\pm 1,3,2))$ is infinite in all cases except the Poincar\'e homology sphere
$S^3/I^*=L(5,3,2).$ For example, $L(7,3,2)$ is the quotient of the universal
cover $\widetilde{SL}(2,\bbr)$ of $SL(2,\bbr)$ by a co-compact
discrete subgroup $\grG\subset \widetilde{SL}(2,\bbr),$ and there is a
covering of $L(7,3,2)$ by the total space of a nontrivial circle
bundle over a Riemann surface of genus $g=3$ \cite{Mil75}. More generally, except for the
Poincar\'e homology sphere, $L(6k\pm 1,3,2)$ has a finite covering by a manifold
that is diffeomorphic to a nontrivial circle bundle over a Riemann
surface of some genus $g>1.$
\end{exmp}

We end this section with a brief discussion of how one can
distinguish integral homology spheres in three dimensions. This
can be done by the Casson invariant $\grl$ which roughly speaking
counts the number of irreducible representations of $\pi_1$ in
$SU(2).$  More explicitly for a Brieskorn homology sphere
$L(p,q,r)$ Fintushel and Stern \cite{FiSt90} proved that $\grl(L(p,q,r))$ equals
$-\frac{1}{2}$ times the number of conjugacy classes of
irreducible representations of $\pi_1(L(p,q,r))$ into $SU(2).$
Furthermore, they showed how $\grl(L(p,q,r))$ can be computed from
the Hirzebruch signature of the parallelizable manifold $V(p,q,r)$
whose boundary is  the link $L(p,q,r).$ Explicitly,
$$\grl(L(p,q,r)) = \frac{{\rm sig}~V(p,q,r)}{8}.$$
Now Brieskorn \cite{Bri66} had computed the signature of $V(6k-1,3,2)$ to be $-8k,$
giving $\grl(L(6k-1,3,2)) = -k.$ Thus, $L(6k-1,3,2)$ give an infinite sequence of
inequivalent homology spheres. Moreover, it is known that $\grl(L(7,3,2))=-1$, and since
$k=1$ from the previous sequence is the Poincar\'e homology sphere, which is clearly
distinct, we see that $L(7,3,2)$ is also distinct from the others.
For further discussion and references we refer the reader to the recent book \cite{Sav02}.

\section{$\eta$-Einstein Metrics in Dimension 5 and Higher}

In this section we give a large number of examples of manifolds in dimension five and
higher that admit Sasakian $\eta$-Einstein metrics. For null and negative Sasakian
structures the proof of existence amounts to presenting examples as links of weighted
homogeneous polynomials and then invoking Theorem \ref{nullnegCal}.
In the positive or Fano case the existence of Sasakian $\eta$-Einstein metrics is equivalent
to the existence of Sasakian-Einstein metrics which has been studied quite
extensively \cite{BGN03a, BGK03, BGKT03, Kol04}.

The list of
spin 5-manifolds which admit many inequivalent families of
Sasakian-Einstein metrics is now long and keeps growing. For
example, $S^5$ has at least 68 inequivalent Sasakian-Einstein
metrics realized as Brieskorn-Pham links. The example below gives a
subfamily.

\begin{exmp}  Consider sequences of the form
$\bfa=(2,3,7,m)$. By explicit calculation, the corresponding link
$L(\bfa)$ gives a Sasakian-Einstein metric on $S^5$ if $5\leq m\leq
41$ and $m$ is relatively prime to at least two of $2,3,7$. This
is satisfied in $27$ cases.
\end{exmp}

Recently, Koll\'ar has shown that there are many inequivalent
Sasakian-Einstein metrics on all $k$-fold connected sums of
$S^2\times S^3$ \cite{Kol04}. On the other hand Boyer and Galicki
showed that \cite{BG03p}

\begin{thm} The links $L(3,3,3,k)$ $k>2$, $L(2,4,4,p)$, $p>2$ and $L(2,3,6,m)$, $m>4$ all
admit
Sasakian-Einstein structures.
\end{thm}
The above links are rational homology 5-spheres under suitable
assumptions on $k,p,m$. For example $L(3,3,3,k)$ is a rational
homology sphere with $H_2(L,\bbz)=\bbz_k\oplus\bbz_k$ as long as
$k$ is prime to 3. Since any $L(3,3,3,k)$ is spin, by Smale's
theorem \cite{Sm62}, the second homology group completely
determines diffeomorphism type of $L(3,3,3,k)$.

In the null case there appear to be only a finite number of manifolds admitting such
structures. We begin with a non-existence result in dimension 5.

\begin{thm}\label{nonull} Let $\cals$ be a null Sasakian structure on a compact 5-manifold
with $H_1(M,\bbz)=0.$ Then $2\leq b_2(M)\leq 21$ and $H_2(M,\bbz)$ is torsion free.
Furthermore, if $b_2(M)=21$ then $\cals$ is regular and $M$ is diffeomorphic to $\#
21(S^2\times S^3)$, and $M/\calf_\xi$ is a K3 surface.
\end{thm}

\begin{proof} The lower bound on $b_2(M)$ is immediate from Corollary 1.10 of
\cite{BGN03a}, while the upper bound as well as the last statement follow from Corollary 81
of \cite{Kol04b}. That $H_2(M,\bbz)$ is torsion free follows from Proposition 80 of
\cite{Kol04b}.
\end{proof}

We mention that the lower bound in Theorem \ref{nonull} holds with the weaker condition
$H_1(M,\bbr)=0.$ We now have

\begin{cor}
$S^5$ and $S^2\times S^3$ or any quotient by a finite group do not admit null Sasakian
structures.
\end{cor}
Compact 5-manifolds with null Sasakian structures can be easily obtained from the list of
95 orbifold K3 surfaces.

\begin{exmp}\label{reid} In 1979 Reid \cite{Rei79} produced a list of 95 weighted K3
surfaces given as well-formed hypersurfaces in weighted projective space
$\bbp(w_0,w_1,w_2,w_3)$. His result generalizes the standard
construction of the K3 surface as a quartic in $\bbc\bbp(3)$. With
each example of Reid one can consider the associated 5-dimensional
link $L_f$. By a result of Boyer and Galicki \cite{BG01b} all  $L_f$
are diffeomorphic to some $k$-fold connected sum of $S^2\times
S^3$, where $k=b_2(L_f)$ can be computed in each case by the Milnor-Orlik procedure
\cite{MiOr70}. For
instance $L(4,4,4,4)\simeq\#21(S^2\times S^3)$ is regular and it
is a circle bundle over the $K3$-surface $X_4\subset\bbp^3$. No.3
On Reid's list in the table of \cite{Fle00} is also a BP link
$L(6,6,6,2)\simeq\#21(S^2\times S^3)$. The degree 6 surface in
$X_6\subset\bbp(1,1,1,3)$ given by vanishing of
$f(\bfz)=z_0^6+z_1^6+z_2^6+z_3^2$ is easily seen to be
non-singular and as a result we also get a regular null Sasakian
structure on a K3-surface and regular null Sasakian structure on
the corresponding circle bundle. However, $X_4$ and $X_6$ cannot
be the same as complex algebraic varieties and
$L(4,4,4,4)$ and $L(6,6,6,2)$ are not equivalent as Sasakian
manifold. In all, there are 11 more examples with are
Brieskorn-Pham links. All Reid's examples have $3\leq b_2(L_f)\leq
21$. Hence, for instance, the link $L_f$ with
$$f(\bfz)=z_0^2+z_1^2+z_2z_0z_1z_3+z_2^2z_3+z_2^3z_0+z_3^3$$
must be diffeomorphic to $\#3(S^2\times S^3)$. No. 14 on the list,
is a BP link $L(2,3,12,12)$ and it diffeomorphic to
$\#20(S^2\times S^3)$. It is interesting that $b_2(L)=17$ does not
occur on the Reid's list. But other examples can be considered.
For instance, Fletcher gives a list of 84 examples of codimension
2 complete intersections \cite{Fle00}. On the other hand, the only
codimension 3 weighted complete intersection is the intersection
of 3 quadrics in $\bbp^5$.
\end{exmp}

By writing a simple Maple program we have computed the second Betti number of all 95
orbifold K3 surfaces, and we find
\begin{cor} $\#k(S^2\times S^3)$ admit both positive and null Sasakian
$\eta$-Einstein structures for $3\leq k\leq21$ and $k\not=17$.
\end{cor}

The question which of the simply connected compact
spin 5-manifolds of Smale's list admit null Sasakian structures is
equivalent to classifying all orbifold K3-surfaces. This is still
open. An interesting question is whether $k=2$ or $17$ can occur. In higher dimensions
there are many, though probably finite, null Sasakian structures, and they will all admit
$\eta$-Einstein null structures. Interesting examples include the over 6000 Calabi-Yau
orbifolds in complex dimension 3 \cite{CLS90}.

In the negative case even just Brieskorn-Pham links produce
already very many examples. The simplest ones are the Fermat hypersurfaces of degree
$d$  in $\bbp^3.$ For $d=4$ it gives $\#21(S^2\times S^3)$ as discussed in Example
\ref{reid} above. For $d\geq 5$ we get negative Sasakian $\eta$-Einstein metrics on
$\#k(S^2\times S^3)$ for $k=(d-2)(d^2-2d+2)+1.$ These begin with $k=52$ and grow
rapidly. To obtain negative Sasakian $\eta$-Einstein structures we turn to some different
examples. First, such structures exist on rational homology spheres, and in particular on
$S^5.$

\begin{exmp}\label{kkp} Consider $L(k,k,k+1,p)$. Such a link is negative as long as
$k,p>3$ or with $k=3$ and $p>12$. Moreover, $L(k,k,k+1,p)$ is a 5-sphere as
long as $p$ is prime to both $k$ and $k+1$ which can easily be
arranged.
\end{exmp}

\begin{prop} Every negative Sasakian structure on $S^5$ is non-regular and there exists
infinitely many inequivalent negative Sasakian
$\eta$-Einstein structures. Hence, $S^5$ admits infinitely many
different Lorentzian Sasakian-Einstein structures.
\end{prop}

\begin{proof}
The existence of infinitely many non-regular negative Sasakian structures on $S^5$ follows
immediately from Example \ref{kkp}. These lie in different deformation classes, and from
Theorem \ref{aub-yau} each deformation class has a negative $\eta$-Einstein metric.
To prove the first statement we assume that $S^5$ has a regular negative Sasakian
structure. Then we have a circle bundle $S^5\ra{1.3} X$ with $c_1(X)$ negative,
$\pi_1(X)=\{1\},$ and $b_2(X)=1.$ This is impossible by a Theorem of Yau, cf. Theorem
V.1.1 of \cite{BPV}.
\end{proof}

It is rather easy to produce examples of negative Sasakian
$\eta$-Einstein structures on $\#k(S^2\times S^3)$; however, we haven't been able to find
one nice sequence which does this, or even a sequence that works for odd $k.$ The
examples below give the essential ideas.

\begin{exmp}\label{kk+1} This is example \ref{kkp} with $p=k+1.$ Now the BP link is
$L(k,k,k+1,k+1)$ and it is easy to see that the induced Sasakian structures are negative as
long as $k\geq 4.$ Here the Milnor-Orlik method gives $b_2(L(k,k,k+1,k+1))=k(k-1).$ This
gives negative Sasakian  $\eta$-Einstein structures on $\#k(S^2\times S^3)$ for infinitely
many $k,$ but with larger and larger gaps as $k$ grows.  For low values of $k,$ for
example, $k=4$ and $5,$ we get negative Sasakian structures on $\#12(S^2\times S^3),$
and $\#20(S^2\times S^3)$, respectively.
\end{exmp}

\begin{exmp}\label{ex2} Consider the BP links $L(p,q,r,pqr)$ such that
$p,q,r\in\bbz^+$ are pairwise relatively prime. The transverse
space is a surface $X_{pqr}\subset\bbp(pr,qr,pq,1)$ and it is
automatically well-formed, i.e., it has only isolated orbifold
singularities. Hence, $L(p,q,r,pqr)$ is diffeomorphic to
$\#k(S^2\times S^3)$, where
$$k=b_2(L(p,q,r,pqr))= (pqr -pq-pr-qr-1) +p+q+r.$$
Now, $L(p,q,r,pqr)$ is negative when the term in parenthesis is positive, i.e. when
$pq+pr+qr+1<pqr$.
To simplify we can consider $L(2,3,r,6r)$ with $r$ prime to both 3
and 2. $L(2,3,7,42)$ is one of the examples on Reid's list. But
from $r>7$ we get infinitely many examples with
second Betti number $b_2=2(r-1).$ For example for $r=11,$  we see that $L(2,3,11,66)$ is
diffeomorphic to $\#20(S^2\times S^3)$.
\end{exmp}

Notice that both Examples \ref{kk+1} and \ref{ex2} produce $\eta$-Einstein metrics on
$\#k(S^2\times S^3)$ for even $k$ only. We haven't yet found a nice series that gives
metrics for odd values of $k.$ However, we can obtain $\eta$-Einstein metrics for odd $k$
if we consider the more general weighted homogeneous polynomials.

\begin{exmp}\label{whpex}
Here we give a few examples only to illustrate the method. A much more extensive list can
be generated with the aid of a computer a la \cite{jk1,BGN03c}.
We begin with two inequivalent `twin' Sasakian $\eta$-Einstein metrics on $S^2\times
S^3.$ Consider the weighted homogeneous polynomials
$$z_0^{21}z_1+z_1^5z_2+z_2^3z_0+z_3^2,\qquad
z_0^{21}z_1+z_1^5z_0+z_2^3z_1+z_3^2.$$
Both have degree $316$ and both give links diffeomorphic to $S^2\times S^3.$ The weight
vectors are $\bfw=(13,43,101,158)$ and $\bfw=(11,61,85,158),$ respectively. So
$|\bfw|-d=-1$ giving negative Sasakian $\eta$-Einstein metrics on $S^2\times S^3.$

Similarly, the polynomials
$$z_0^{20} +z_1^3z_3+z_2^3z_1+z_3^2z_0\qquad
z_0^{10}+z_1^4z_2+z_2^2z_3+z_3^2z_0$$
have degree $40$ and $159$, respectively. They
give negative Sasakian $\eta$-Einstein metrics on $\#7(S^2\times S^3)$ and
$\#2(S^2\times S^3).$
\end{exmp}

We have

\begin{cor} $\#7(S^2\times S^3)$, $\#12(S^2\times S^3)$ and $\#20(S^2\times S^3)$ all
admit positive, null, and
negative Sasakian $\eta$-Einstein structures.
\end{cor}

It is easy to construct examples of negative Sasakian
$\eta$-Einstein metrics on homotopy spheres in arbitrary odd
dimensions.

\begin{exmp} Let the integers $r_i$, $i=1,\ldots,2m$ be any pairwise prime positive integers.
Then
the $(4m+1)$-dimensional link $L(2,2r_1,\ldots,2r_{2m},a)$ is
diffeomorphic to the standard sphere if $a\equiv \pm 1 \mod 8$ and
to the Kervaire sphere if $a\equiv \pm 3 \mod 8.$ The link will be
negative when $\sum_{i}\frac{1}{r_i}<\frac{a-2}{a}$ which is
easily satisfied for $r_i$'s large enough.
\end{exmp}

\begin{cor} In any dimension $4m+1$ both the standard and the
Kervaire spheres admit infinitely many inequivalent negative
Sasakian $\eta$-Einstein structures.
\end{cor}

In dimension $4m+3$ one can construct examples of negative
Sasakian structures should exists on all homotopy spheres which
bound parallelizable manifolds. That is to say, all BP links that
by Brieskorn Graph Theorem are homeomorphic to a sphere of
dimension $4m+3$ and are negative should easily contain all
possible oriented diffeomorphism types. We have checked it in
dimension 7 with the following

\begin{thm}
All 28 oriented diffeomorphism types of
homotopy 7-spheres admit negative Sasakian structures. Hence, they all admit Lorentzian
Sasakian-Einstein structures.
\end{thm}

\begin{proof}
Consider $L(k,k,k,k+1,p)$, where $p$ is prime to both
$k$ and $k+1$. The link $L$ is a homotopy 7-sphere and it is
negative for $k$ and $p$ large enough. Using the computer codes of \cite{BGKT03}, it is
easy to check that such links realize all 28 oriented diffeomorphism types.
\end{proof}

Similarly, one should be able to show that the homotopy spheres in
dimensions 11 and 15 that bound parallelizable manifolds admit negative Sasakian structures.

\begin{rem} Many of the Sasakian $\eta$-Einstein structures admit large moduli spaces.
This is particularly true for BP type polynomials. We refer the reader to our previous work
\cite{BGN03c,BGK03,BGKT03}.
\end{rem}

\bigskip
\section{Positive $\eta$-Einstein Structures and Einstein-Weyl Geometry}
\medskip

In this chapter we investigate relations between Einstein-Weyl
structures and $\eta$-Einstein metrics. The first insight into this
relation was provided by Swann and Pedersen \cite{PeSw93} and by Higa \cite{Hig93}
who studied Einstein-Weyl geometries of circle bundles over positive
K\"ahler-Einstein manifolds.  Later Narita in a series of papers
made the connection to $\eta$-Einstein geometry more explicit \cite{Nar93,Nar97,Nar98}.
We begin with some basic
facts about Einstein-Weyl geometry (see \cite{CalPe99} for a
review of the subject).

\begin{defn} A {\it Weyl structure} on a manifold $M$ of dimension $n\geq3$
is defined by a pair $\calw=([g],D)$, where $[g]$ is a conformal
class of Riemannian metrics and $D$ is the unique torsion-free connection preserving
$[g]$. The connection $D$ is called the {\it Weyl connection}.
\end{defn}

Then, for every Riemannian metric $g$ in the conformal class,
there exists a $1$-form $\theta$, uniquely determined by $D$ and
$g$, such that, for every  $X,Y \in \CX(M)$,
\begin{equation}\label{connection}
D_{X}Y=\nabla_{X}Y+\theta(X)Y+\theta(Y)X-g(X,Y)\zeta,
\end{equation}
where $\nabla$ denotes the Levi-Civita connection of $g$ and
$\theta^\#$ the dual vector field of $\theta$ with respect to $g$. In
particular we have
\begin{equation}\label{fondam}
Dg=-2\theta\otimes g.
\end{equation}
The above condition is invariant under Weyl transformation, i.e.,
\begin{equation}\label{Weyl}
g'=e^{2f}g,\qquad \theta'=\theta+df,\qquad f\in C^\infty(M)
\end{equation}

Hence, choosing a Riemannian metric in $[g]$ we sometimes abuse the terminology and
refer to the pair $\calw=(g,\theta)$ as a Weyl structure. It is
known that on a compact, conformal manifold of dimension at least
$3$, the conformal class $[g]$ contains a unique (up to homothety)
metric $g$ such that the corresponding form $\theta$ is
$g$-co-closed (cf. \cite{Gau95}). This metric is {\it the Gauduchon
metric} or {\it the Gauduchon gauge}.

\begin{defn} A conformal manifold
$(M,[g],D)$ is called {\it Einstein-Weyl} if
$${\rm Ric}^D(X,Y)+{\rm Ric}^D(Y,X)=\Lambda g(X,Y)$$
for some smooth function $\Lambda$ on $M$. \end{defn}

Clearly the Einstein-Weyl condition on the symmetrized Ricci
tensor of $D$ is conformally invariant. A Weyl structure is said
to be {\it closed} (resp. {\it exact}) if the Weyl connection is
locally (resp. globally) the Levi-Civita connection of a
compatible metric. Hence a  closed Einstein-Weyl structure admits
local (but not necessarily global) compatible Einstein metrics.

On compact manifolds of dimension at least $3$ with a closed but
not exact Einstein-Weyl structure the form $\theta$ of the
Gauduchon metric $g$ is $g$ parallel: $\nabla\theta=0$
\cite{Gau95}. In particular, $\theta$ is closed, hence
$g$-harmonic.

On an almost contact metric manifold $(M,\xi,\eta,\Phi, g)$  it is
natural to ask whether $M$ has an Einstein-Weyl structure
$\calw=(g,\theta)$ such that $\theta=f\eta$, for some function
$f\in C^\infty(M)$, and whether the metric $g$ compatible with
the almost contact structure is the Gauduchon metric of
$[g]$. In fact, if $(M,\xi,\eta,\Phi,g)$ is an almost contact manifold
then, for any 1-form $\theta$, the Weyl structure
$\calw=(g,\theta)$ on $M$ is connected to
$\cals=(\xi,\eta,\Phi,g)$ by several fundamental relations. The
existence of Einstein-Weyl structures on almost contact metric
manifolds have been considered by several authors \cite{Matz02},
\cite{Matz00}, \cite{Nar97}, \cite{Nar98}. Here we review and then
extend some of the known results.

Let $(M,\xi,\eta,\Phi,g)$ be a $(2n+1)$-dimensional, $n\geq 1$,
almost contact metric manifold. If $\calw=(g,\theta)$ is a Weyl
structure on $M$, for every  $X,Y \in \CX(M)$, the Ricci tensor
${\rm Ric}^D$ of $D$ and ${\rm Ric}_g$ of $\nabla$ are related by
the following equation \cite{Hig93}:
\begin{equation}\label{Ricci}
\begin{array}[t]{l}
{\rm Ric}^D(X,Y)={\rm
Ric}_g(X,Y)-2n(\nabla_{X}\theta)(Y)+(\nabla_{Y}\theta)(X)+
\\ \ \\
+(2n-1)\theta(X)\theta(Y)+(\delta
\theta-(2n-1)|\theta|^{2})g(X,Y),
\end{array}
\end{equation}
where $\delta \theta$,  $|\theta|$ are the codifferential  and the
point-wise norm of $\theta$ with respect to $g$.
The following local characterizations of Einstein-Weyl structures
can be found in \cite{Hig93}.
\begin{prop}\label{EW1}
Let $\calw=(g,\theta)$ be a Weyl structure on an almost contact
metric manifold $(M,\xi,\eta,\Phi,g)$ of dimension $2n+1\geq 3$.
Then  $\calw=(g,\theta)$ is an Einstein-Weyl structure if and only
if there exists a smooth function $\sigma$ on $M$ such that:
\begin{equation}\label{prop}
\frac{1-2n}{2}((\nabla_{X}\theta)(Y)+(\nabla_{Y}\theta)(X))+(2n-1)\theta(X)\theta(Y)
+{\rm Ric}_g(X,Y) = \sigma g(X,Y),
\end{equation}
for every $ X,Y \in \CX(M)$.
\end{prop}

\begin{prop}\label{EW2}
Let $(M,\xi,\eta,\Phi,g)$ be an almost contact metric manifold of
dimension $2n+1\geq 3$ and $\theta$ a 1-form on $M$. Consider the
Weyl structures $\calw^{\pm}$ defined by $\calw^{\pm}=(g,\pm
\theta)$. Then both $\calw^{+}$ and $\calw^{-}$ are Einstein-Weyl
if and only if $(g,\theta)$ satisfies the following two equations
for every $X,Y \in \CX(M)$:
\begin{equation}\label{prima}
(\nabla_{X}\theta)(Y)+(\nabla_{Y}\theta)(X)+\frac{2}{2n+1}\delta
\theta \ g(X,Y)=0,
\end{equation}

\begin{equation}\label{seconda}
{\rm
Ric}_g(X,Y)-\frac{s}{2n+1}g(X,Y)=\frac{2n-1}{2n+1}|\theta|^{2}\
g(X,Y)-(2n-1) \theta(X) \theta(Y),
\end{equation}
where $s$ is the scalar curvature of $g$.
\end{prop}

First of all, as a direct consequence of (\ref{prop}), we can
state

\begin{prop}\label{trivial} Let $(M,\xi,\eta,\Phi, g)$ be a
$K$-contact metric manifold of dimension $2n+1\geq 3$ and suppose
that $M$ admits an Einstein-Weyl structure $W=(g, \theta)$ with
$\theta=f \eta$ for some non constant function $f$ on $M$. Then
$X(f)=0$ for every vector field $X$ on $M$ such that $\eta(X)=0$
and the Ricci tensor of $M$ is given by ${\rm Ric}_g=\sigma
g-(2n-1)(f^{2}-\xi(f))\eta \otimes \eta$.
\end{prop}

\begin{proof}
We recall, at first, that the Ricci tensor of a $K$-contact
manifold always satisfies the relation ${\rm
Ric}_g(X,\xi)=2n\eta(X)$. Then, for every $X$ horizontal vector
field  on $M$ (here and in the following we shall call {\it
horizontal} any vector field on $M$ orthogonal to $\xi$) and
$Y=\xi$, (\ref{prop}) becomes
\begin{equation}\label{easy}
(\nabla_{X}f \eta)(\xi)+(\nabla_{\xi}f \eta)(X)=0,
\end{equation}
which easily implies the first part of the proposition. A
straightforward computation proves the second part of the assert.
\end{proof}

\begin{cor}
If a $K$-contact manifold $M$ of dimension at least 5 admits an
Einstein-Weyl structure $\calw=(g, f\eta)$ with $f\in
C^\infty(M)$, then $M$ is $\eta$-Einstein. In particular, in such
case, $\sigma$ and $f^{2}-\xi(f)$ must be constant.
\end{cor}

The following example generalizes the 3-dimensional case of
Example 35.

\begin{exmp} \rm  Let $H(n)$ be the Heisenberg Lie group
\begin{equation} {H(n)=\left\{\left(
\begin{array}{ccc}
1 &\bfx^t & z\\
0 & 1 & \bfy\\
0 & 0 & 1 \end{array}\right)\  ;\
\bfx,\bfy\in\bbr^n,z\in\bbr\right\}},
\end{equation}
and let $g$ be the following left invariant metric on $H(n)$
\begin{equation}\label{metr2}
g=d\bfx\cdot d\bfx+d\bfy\cdot d\bfy+(dz-\bfx\cdot d\bfy)^2.
\end{equation}

The Sasakian structure
$\cals=(\xi,\eta,\Phi,g)$ on $H(n)$ is defined by the formulas
\begin{equation} \left\{
\begin{array}{lcl}
\eta= dz-\bfx\cdot d\bfy,\ \  \xi =\frac{\partial}{\partial z}  \\
\Phi =\sum_i\bigl((\frac{\partial}{\partial x^i}+y^i\frac{\partial}{\partial z})\otimes dy^i
-\frac{\partial}{\partial y^i}\otimes dx^i\bigr).
\end{array}\right.
\end{equation}
One can easily check that this is a null Sasakian $\eta$-Einstein
structure on $H(n)$ \cite{Sas65}. Moreover, a transverse homothety of $\cals$ by
a constant $a$ gives a null Sasakian $\eta$-Einstein structure
$\bar{\cals}=(\bar{\xi},\bar{\eta},\bar{\Phi},\bar{g})$ which admits the
Einstein-Weyl structure $\calw=(\bar{g},f \bar{\eta})$, with
$f=\alpha\tan(z+c), \alpha,c\in \bbr$, provided that
$a=\frac{1}{\alpha}$, and $\alpha>0$ such that
$\alpha^2=\frac{2n+2}{2n-1}$. This generalizes to arbitrary odd dimension the
Einstein-Weyl structure described in \cite{Nar97,Nar98}. In addition, one can get compact
examples by taking quotients of $H(n)$ by discrete subgroups.
\end{exmp}

From the
formula (\ref{prop}) we can also deduce that, if $\xi$ is a Killing
vector field on the almost contact metric manifold
$(M,\xi,\eta,\Phi, g)$ and $M$ admits an Einstein-Weyl structure
$\calw^{+}=(g,\mu \eta)$ for some constant $\mu$, then $M$ admits
the Einstein-Weyl structure $\calw^{-}=(g,-\mu \eta)$ too. More
precisely, applying the Proposition \ref{EW2} to the $K$-contact
manifolds, we obtain the following

\begin{thm}\label{Kc} Let $(M,\xi,\eta,\Phi, g)$ be a
$K$-contact manifold of dimension $2n+1\geq 3$ and $\theta$ a
1-form on $M$. Suppose that both $\calw^{\pm}=(g,\pm \theta)$ are
Einstein-Weyl structures on $M$. Then $\theta$ is co-closed, but not closed, so $g$ is the
Gauduchon metric, the dual vector field $\theta^\#$ to $\theta$ is a Killing vector field, and
either
\begin{enumerate} \item $\theta(\xi)=0$ or
\item $\theta=\mu \eta$ with $\mu \in\bbr$.
\end{enumerate}
\end{thm}

\begin{proof}
In fact, if both $\calw^{\pm}=(g,\pm \theta)$ are Einstein-Weyl
structures on $(M,\xi,\eta,\Phi, g)$, the 1-form $\theta$
satisfies the both equations of Proposition \ref{EW2} . In particular, substituting
$\xi$ for $Y$ in (\ref{seconda}) we get
\begin{equation}\label{nuova}
2n\eta(X)-\frac{s}{2n+1}\eta(X)=\frac{2n-1}{2n+1}|\theta|^{2}\
\eta(X)-(2n-1) \theta(X) \theta(\xi).
\end{equation}
Then the relation
\begin{equation}\label{fond}
\theta(X) \theta(\xi)=0
\end{equation} is true for every horizontal vector field $X$.
Now define the set
$$S=\{x\in  M~|~\theta(X)_x=0~\forall~{\rm horizontal}~X\}.$$
$S$ is closed, and by (\ref{fond}) $\theta(\xi)=0$ on the open set $M\backslash S.$
Then putting $X=Y=\xi$ in (\ref{prima}) implies that
$\delta \theta=0$ on $M\backslash S.$ But then again (\ref{prima}) implies
$(\nabla_{X}\theta)(Y)+(\nabla_{Y}\theta)(X)=0$ on $M\backslash S$ for any $X,Y \in
\CX(M)$. This means that the dual vector field $\theta^{\#}$ is a Killing field on
$M\backslash S.$ Consequently, if $\theta$ were closed on $M\backslash S$, it would also
be parallel with respect to the Levi-Civita connection $\nabla$ of $g$. In
particular for every $X \in \CX(M)$ we have
\begin{equation}\label{nulla}
0=(\nabla_{X}\theta)(\xi)= -\theta(\Phi (X))
\end{equation}
implying that $\theta$ identically vanishes on $M\backslash S.$
This proves that $\theta$ cannot be closed on $M\backslash S.$

On the other hand, the set $S_0$ where $\theta(\xi)\neq 0$ is open and a subset of $S$
by (\ref{fond}). Thus,  $\theta=\mu \eta$ on $S$ for some $\mu\in
C^\infty(S)$. So, for $X,Y$ horizontal vector fields
(\ref{prima}) becomes
\begin{equation}\label{terza}
\mu ((\nabla_{X}\eta)(Y)+
(\nabla_{Y}\eta)(X))+\frac{2}{2n+1}\delta \theta \ g(X,Y)=0
\end{equation}
on $S.$ But the first term vanishes since $\xi$ is a Killing vector field. This implies that
$\delta \theta=0$ on $S$ and hence, on all of $M.$ But then (\ref{prima}) gives
$(\nabla_{X}\theta)(Y)+(\nabla_{Y}\theta)(X)=0,$ on all $M$ for all $X,Y \in
\CX(M)$ which means that $\theta^\#$ is a Killing vector field, and also implies that $\mu$
is a constant. This completes the
proof of the theorem.
\end{proof}

We now discuss several corollaries of Theorem \ref{Kc}.
First combining our results with \cite{Gau95} we have

\begin{cor} \label{gau}
Let $(M,\xi,\eta,\Phi, g)$ be a compact $K$-contact manifold of
dimension $2n+1\geq 3$ and $\theta$ a 1-form on $M$. Given the
Weyl structure $\calw=(g,\theta)$, then the metric $g$ is the
Gauduchon metric of the conformal class $[g]$ if and only if $M$
admits both $\calw^{\pm}=(g,\pm \theta)$ as Einstein-Weyl
structures.
\end{cor}

\begin{cor} \label{EWKc} Let $(M,\xi,\eta,\Phi,g)$ be a
$K$-contact manifold of dimension $2n+1\geq 5$. Then
$(M,\xi,\eta,\Phi, g)$ admits an Einstein-Weyl structure
$\calw=(g,\theta)$ with $\theta=\mu \eta$, $\mu \in \bbr$, if and
only if $M$ is $\eta$-Einstein with Einstein constants
$\lambda,\nu$ such that $\nu<0$.
\end{cor}
\begin{proof}
If $(M,\xi,\eta,\Phi,g)$ admits an Einstein-Weyl
structure with 1-form $\theta=\mu \eta$ for some constant $\mu$, then
from (\ref{prop}) we obtain
\begin{equation}\label{prop1}
{\rm Ric}_g =\sigma g-(2n-1)\mu^{2} \eta \otimes \eta ,
\end{equation}
which proves the only if part of the assertion. On the contrary, if we
suppose that the Ricci tensor of $M$ is given by
\begin{equation}\label{rho}
{\rm Ric}_g=\lambda g+\nu \eta \otimes \eta,
\end{equation}
where $\nu<0$, (\ref{prop}) is satisfied for $\sigma= \lambda$ and
$\theta=\mu \eta$ where $\mu$ is a constant such that
$\mu^{2}=-\frac{\nu}{2n-1}$.
\end{proof}

In the case of Sasakian structures Theorem \ref{Kc} can be improved:

\begin{thm}\label{Sasakian} A Sasakian manifold $(M,\xi,\eta,\Phi,g)$
of dimension $2n+1\geq 5$ admits both the Einstein-Weyl structures
$\calw^{+}=(g,\theta)$ and $\calw^{-}=(g,-\theta)$ for some 1-form
$\theta$ if and only if is $\eta$-Einstein with Einstein constants
$(\lambda,\nu)$ such that $\nu<0$.
\end{thm}
\begin{proof}
From the proof of Theorem \ref{Kc} we know that, if
$(M,\xi,\eta,\Phi,g)$ admits both $\calw^{\pm}$ as Einstein-Weyl
structures, then for the 1-form $\theta$ the relation $\theta(X)
\theta(\xi)=0$ holds, for every horizontal vector field $X$ on
$M$. On the other hand, since on a Sasakian manifold the Ricci
tensor satisfies the equality ${\rm Ric}_g(\Phi(X),\Phi(Y))={\rm
Ric}_g(X,Y)-2n\eta(X) \eta(Y)$, substituting in (\ref{seconda})
$\Phi(X)$ for $X$ and $\Phi(Y)$ for $Y$ with $X,Y$ horizontal
vector fields, we get
\begin{equation}\label{uguale}
{\rm
Ric}_g(X,Y)-\frac{s}{2n+1}g(X,Y)=\frac{2n-1}{2n+1}|\theta|^{2}\
g(X,Y)-(2n-1) \theta(\Phi(X)) \theta(\Phi(Y)),
\end{equation}
where $s$ is the scalar curvature of $g$. The comparison between
(\ref{seconda}) and this last equation implies that $\theta$ must also
obey the equation
\begin{equation}\label{ultima}
\theta(\Phi(X)) \theta(\Phi(Y))=\theta(X) \theta(Y),
\end{equation}
for all horizontal vector fields $X,Y$ on $M$. Now, arguing as in the proof of Theorem
\ref{Kc}, if we suppose
$\theta(\xi)=\eta(\theta^\#)=0$ and consider $X=Y=\theta^\#$ in (\ref{ultima}), we
obtain $|\theta|=0$ so that $\theta$ identically vanishes on $M$.
Finally, the Corollary \ref{EWKc} implies the theorem.
\end{proof}

\begin{rem} We remark that, by taking Corollary \ref{gau} into account,
in the case of a compact Sasakian manifold $(M,\xi,\eta,\Phi,g),$
the above Theorem \ref{Sasakian} provides the necessary and
sufficient condition for $g$ to be the Gauduchon metric of
the conformal class $[g]$ for some Einstein-Weyl structure on $M$.
\end{rem}

Finally, the following corollary follows easily from the previous analysis.

\begin{cor}\label{squashed3}
Let $(\xi,\eta,\Phi,g)$ be a Sasakian-Einstein structure on a manifold $M.$ Then $M$
admits a pair of Einstein-Weyl structures
$\calw^{\pm}=(g,\pm \mu\eta)$ obtain by squashing the \Se structure to an
$\eta$-Einstein structure
with constants $(\lambda,\nu)$ and $\nu<0.$ In this case
$\mu^2=-\frac{\nu}{2n-1}$. In particular, all homotopy spheres of dimension $4n+1,7,11$
and $15$ that bound parallelizable manifolds admit many pairs of Einstein-Weyl structures.
\end{cor}

The above corollary combined with recent results of \cite{BGN03c, BGK03, BGKT03,
BG03p, Kol04, Kol04b}
establishes the existence of Einstein-Weyl structures on many other compact
simply connected spin manifolds in odd dimensions.


\def\cprime{$'$} \def\cprime{$'$} \def\cprime{$'$}
\providecommand{\bysame}{\leavevmode\hbox to3em{\hrulefill}\thinspace}
\providecommand{\MR}{\relax\ifhmode\unskip\space\fi MR }
\providecommand{\MRhref}[2]{%
  \href{http://www.ams.org/mathscinet-getitem?mr=#1}{#2}
}
\providecommand{\href}[2]{#2}

\end{document}